\documentclass[11pt,oneside]{book}
\usepackage{geometry}                		
\usepackage{graphicx}									
\usepackage{amssymb}\vspace{0.25cm}
\usepackage{mathtools}
\usepackage{setspace}
\usepackage{tikz-cd}
\usepackage[mathscr]{euscript}
\newcommand{\abs}[1]{\left\vert#1\right\vert}
\usepackage[normalem]{ulem}
\usepackage{manfnt}
\usepackage{pgfplots}
\usepackage{pifont}
\usepackage[safe]{tipa}
\usepackage{marvosym}
\usepackage{scalerel,stackengine}
\usepackage{titlesec}
\usepackage{makeidx}
\usepackage{centernot}
\usepackage{xspace}
\usepackage[noauto]{chappg}
\usepackage[OT2,T1]{fontenc}
\usepackage{hyperref}
\usepackage{amsmath}
\usepackage[amsmath,thmmarks]{ntheorem}
\usepackage{tabularx}
\usepackage{yfonts}
\newcolumntype{Y}{>{\raggedright\arraybackslash}X}

\newcommand\mA{%
$A$\xspace
}

\newcommand\mC{%
$C$\xspace
}

\newcommand\mE{%
$E$\xspace
}

\newcommand\mK{%
$K$\xspace
}

\newcommand\mN{%
$N$\xspace
}

\newcommand\mP{%
$P$\xspace
}

\newcommand\mS{%
$S$\xspace
}
\newcommand\mT{%
$T$\xspace
}

\newcommand\mX{%
$X$\xspace
}
\newcommand\mY{%
$Y$\xspace
}
\newcommand\mZ{%
$Z$\xspace
}
\newcommand\pe{%
\mathfrak{p}\xspace
}
\newcommand\A{%
\mathbb{A}
}
\newcommand\C{%
\mathbb{C}
}

\newcommand\M{%
\mathbb{M}
}
\newcommand\N{%
\mathbb{N}
}

\newcommand\Q{%
\mathbb{Q}
}
\newcommand\R{%
\mathbb{R}
}

\newcommand\bx{%
\textbf{x}\xspace
}

\newcommand\elln{%
\ell\text{n}
}

\newcommand\sB{%
\mathcal{B}
}

\newcommand\sE{%
\mathcal{E}
}
\newcommand\sF{%
\mathcal{F}
}

\newcommand\sP{%
\mathcal{P}
}
\newcommand\sS{%
\mathcal{S}
}

\newcommand\alg{%
\text{alg}
}

\newcommand\td{%
\text{d}
}

\newcommand\tI{%
\text{I}
}

\newcommand\tT{%
\text{T}
}

\newcommand\ra{%
\rightarrow
}

\newcommand\ds{%
\displaystyle
}

\newcommand\fP{%
\frak{P}
}
\newcommand\fR{%
\frak{R}
}

\newcommand\modx{%
\text{mod}\hspace{0.05cm}
}
\newcommand\pro{%
\text{pro}%\hspace{0.05cm}
}
\newcommand\vol{%
\text{vol}%\hspace{0.05cm}
}
\newcommand\sgn{%
\text{sgn}\hspace{0.05cm}
}

\newcommand\un[1]{%
\underline{#1}\xspace
}

\newcommand\ov[1]{%
\overline{#1}
}

\newcommand\restr[2]{%
{#1}|{#2}
}

\newcommand\hsx{%
\hspace{0.05cm}
}

\newif\ifcomments
\commentstrue
\newif\ifdmc
\dmcfalse
\newif\ifdmcx
\dmcxtrue
\newif\ifgw
\gwfalse
\usepackage{scalerel,stackengine}
\stackMath
\newcommand\reallywidehat[1]{%
\savestack{\tmpbox}{\stretchto{%
  \scaleto{%
    \scalerel*[\widthof{\ensuremath{#1}}]{\kern-.6pt\bigwedge\kern-.6pt}%
    {\rule[-\textheight/2]{1ex}{\textheight}}%WIDTH-LIMITED BIG WEDGE
  }{\textheight}% 
}{0.5ex}}%
\stackon[1pt]{#1}{\tmpbox}%
}
\makeatletter
\DeclareRobustCommand\widecheck[1]{{\mathpalette\@widecheck{#1}}}
\def\@widecheck#1#2{%
    \setbox\z@\hbox{\m@th$#1#2$}%
    \setbox\tw@\hbox{\m@th$#1%
       \widehat{%
          \vrule\@width\z@\@height\ht\z@
          \vrule\@height\z@\@width\wd\z@}$}%
    \dp\tw@-\ht\z@
    \@tempdima\ht\z@ \advance\@tempdima2\ht\tw@ \divide\@tempdima\thr@@
    \setbox\tw@\hbox{%
       \raise\@tempdima\hbox{\scalebox{1}[-1]{\lower\@tempdima\box
\tw@}}}%
    {\ooalign{\box\tw@ \cr \box\z@}}}
\makeatother

\newtheoremstyle{xx}% name of the style to be used
  {4pt}% measure of space to leave above the theorem. E.g.: 3pt
  {0pt}% measure of space to leave below the theorem. E.g.: 3pt
  {\upshape}% name of font to use in the body of the theorem
  {\bfseries}% name of head font
  {}% punctuation between head and body
  { }% space after theorem head; " " = normal interword space \footnotesize
  {}
  
\makeatletter 
 \newtheoremstyle{myu}%
  {\upshape\item[ \indent\indent\bf\underline{\theorem@headerfont ##2:}]}%
\makeatother
\makeatletter 
 \newtheoremstyle{myn}%
  {\item[\hskip\labelsep \ \bf ##1 \theorem@headerfont ##2.]}%
\makeatother
\theoremstyle{myn}
\newtheorem{theoremn}{Theorem} %<-- Normal Theorem Definition
\theoremstyle{myu}
{\upshape}
\newtheorem{x}[theoremn]{}%<-- Underlined Theorem Definition

 \newtheoremstyle{mr}%
  {\upshape\item[ \indent{\theorem@headerfont ##2. \hspace{.2cm}}]}%
\makeatother
\theoremstyle{mr}
{\upshape}
\title{Periods And Real Numbers}
\author{Garth Warner (Emeritus)\\
Department of Mathematics\\
University of Washington}
\date{}									% Activate to display a given date or no date
\titleformat{\chapter}[display]
{\normalfont\filcenter\huge\bfseries}{}{0pt}{\large}
\titleformat{\chapter}[display]
{\normalfont\filcenter\huge\bfseries}{}{0pt}{\large}
\usepackage[OT2,OT1]{fontenc} 
\newcommand\cyr
{
\renewcommand\rmdefault{wncyr} 
\renewcommand\sfdefault{wncyss} 
\renewcommand\encodingdefault{OT2} 
\normalfont
\selectfont
}

\DeclareTextFontCommand{\textcyr}{\cyr}

\makeindex 
\begin{document}

\maketitle                              % Print title page.

\titlespacing*{\chapter}{0pt}{-50pt}{40pt}
\setlength{\parskip}{0.1em}
\renewcommand{\thepage}{\roman{page}}
\chapter{
ABSTRACT}
\setlength\parindent{2em}
\setcounter{theoremn}{0}
%\renewcommand{\thepage}{\arabic{page}}
%%----------------------------------------------------------------------------------------------01

\ \indent 

The purpose of this book is to provide an introduction to period theory and then to place it within the matrix of recursive function theory.

%%%%%%%%%%%%%%%%%%%%%%%%%%%%%%%%%%%%%%
%%%%%%%%%%%%%%%%%%%%%%%%%%%%%%%%%%%%%%
%%%%%%%%%%%%%%%%%%%%%%%%%%%%%%%%%%%%%%

\chapter{
ACKNOWLEDGEMENT}
\setlength\parindent{2em}
\setcounter{theoremn}{0}
%\renewcommand{\thepage}{\S0-\arabic{page}}
%%----------------------------------------------------------------------------------------------01

\ \indent 

Many thanks to Judith Clare Salzer for typing the manuscript on an IBM Selectric.  
Many thanks to David Clark for converting the typed manuscript to AMS-TeX.

%%%%%%%%%%%%%%%%%%%%%%%%%%%%%%%%%%%%%%
%%%%%%%%%%%%%%%%%%%%%%%%%%%%%%%%%%%%%%
%%%%%%%%%%%%%%%%%%%%%%%%%%%%%%%%%%%%%%

%\include{_acknowledgements}

\centerline{\textbf{\LARGE CONTENTS}}
\vspace{0.5cm}

%%-----------------------------------------------------------------------------------
\allowdisplaybreaks
\begin{align*}
\S1.  \qquad &\text{SEMIALGEBRAIC SETS} 
\\[10pt]
\S2.  \qquad &\text{PERIODS: DEFINITION AND EXAMPLES}  
\\[10pt]
\S3.  \qquad &\text{PERIODS: ALGEBRAIC CONSIDERATIONS} 
\\[10pt]
\S4.  \qquad &\text{PERIODS: THEORETICAL CONSIDERATIONS}
\\[10pt]
\S5.  \qquad &\text{PRIMITIVE RECURSIVE FUNCTIONS}
\\[10pt]
\S6.  \qquad &\text{ELEMENTARY FUNCTIONS}
\\[10pt]
\S7.  \qquad &\text{HIERARCHIES}
\\[10pt]
\S8.  \qquad &\text{COMPUTABILITY}
\\[10pt]
\S9.  \qquad &\text{THE SKORDEV CRITERION}
\\[10pt]
\S10.  \qquad &\text{TECHNICALITIES }
\\[10pt]
\S11. \qquad &\text{NUMERICAL EXAMPLES}
\\[10pt]
\S12.  \qquad &\text{$\R_{\sE^3}$ VERSUS $\R_{\sE^2}$}
\\[10pt]
\S13.  \qquad &\text{RECURSIVE FUNCTIONS}
\\[10pt]
\S14.  \qquad &\text{EXPANSION THEORY}
\\[10pt]
\end{align*}
\vspace{0.5cm}

%%----------------------------------------------------------------------------------------------
%%----------------------------------------------------------------------------------------------

\ \indent 

According to M. Kontsevich and D. Zagier:

\[
\text{* \ * \ * \ * \ * \ * \ *}
\]

``A \un{period} is a complex number whose real and imaginary parts are values of absolutely convergent integrals 
of rational functions with rational coefficients over domains in $\R^n$ given by polynomial inequalities 
with rational coefficients.''
\\
\vspace{11cm}

\un{Ref:} \ 
\textit{Mathematics Unlimited - 2001 and Beyond}, Springer, pp. 771-808.

%%%%%%%%%%%%%%%%%%%%%%%%%%%%%%%%%%%%%%
%%%%%%%%%%%%%%%%%%%%%%%%%%%%%%%%%%%%%%
%%%%%%%%%%%%%%%%%%%%%%%%%%%%%%%%%%%%%%

%\include{__preface}
\pagenumbering{bychapter}
%\chapter{}
%\newpage
%\setcounter{page}{1}
%\renewcommand{\thepage}{1-\arabic{page}}
%\renewcommand{\thepage}{1-\arabic{page}}
\setcounter{chapter}{0}
%%----------------------------------------------------------------------------
\chapter{
$\boldsymbol{\S}$\textbf{1}.\quad  SEMIALGEBRAIC SETS}
\setlength\parindent{2em}
\setcounter{theoremn}{0}
\renewcommand{\thepage}{\S1-\arabic{page}}
%%----------------------------------------------------------------------------------------------01

\begin{x}{\small\bf NOTATION} \ %01
Put
\[
\R_{\alg} \ = \ \R \ \cap \ \ov{\Q}.
\]
Therefore $\R_{\alg}$ is the field of real algebraic numbers, so 
\begin{align*}
\Q \ \subset \ 
&\R_{\alg} \subset \ov{\Q}
\\[1pt]
& \cap 
\hspace{.85cm}
\cap
\\[1pt]
& \hspace{.1cm}
\R \hspace{.3cm} \subset \ \C.
\end{align*}
%\vspace{0.2cm}

[Note: \ 
Viewed as a vector space over $\Q$, $\R_{\alg}$ is infinite dimensional (the algebraic numbers $\sqrt{p}$ ($p$ a prime) are independent over $\Q$).]
\end{x}

\begin{x}{\small\bf \un{N.B.}} \ %02
 $\R_{\alg}$ is a real closed field and  $\R_{\alg}$ is the real closure of $\Q$.  
\end{x}

[Note: \ 
 $\R_{\alg}$ carries the relative topology per $\R$, hence is totally disconnected, i.e., its connected components are points.]
\vspace{0.3cm}

\begin{x}{\small\bf DEFINITION} \ %03
A subset $X \subset \R^n$ is \un{semialgebraic} if there are natural numbers $r$ and $s_1, \ldots, s_r$ and polynomials
\[
 f_{i j}, \hsx  g_{i j} \ \in \ \R_{\alg} [\tT_1, \ldots, \tT_n],
\]
where $1 \leq i \leq r$ and $1 \leq j \leq s_i$,  such that
\[
X 
\ = \ 
\bigcup\limits_{i=1}^r \hsx 
\bigcap\limits_{j=1}^{s_i}  \hsx 
\bigg\{\bx \in \R^n : f_{i j} (\bx) = 0\bigg\} 
\cup 
\bigg\{\bx \in \R^n : g_{i j} (\bx) > 0\bigg\}
.
\]
\end{x}
\vspace{0.2cm}

\begin{x}{\small\bf EXAMPLE} \ %04
$\R^n$ is semialgebraic.
\end{x}
\vspace{0.2cm}

\begin{x}{\small\bf REMARK} \ %05
An algebraic set in $\R^n$ defined using polynomials with coefficients in $\R_{\alg}$ is semialgebraic.
\end{x}
\vspace{0.3cm}
%%----------------------------------------------------------------------------------------------02

\begin{x}{\small\bf LEMMA} \ %06
The semialgebraic sets are closed under the formation of finite unions and finite intersections and are closed under complementation.
\end{x}
\vspace{0.3cm}

\begin{x}{\small\bf LEMMA} \ %07
If $X \subset \R^n$ is semialgebraic and $Y \subset \R^m$ is semialgebraic, then 
$X \times Y \subset \R^n \times \R^m$ is semialgebraic.
\end{x}
\vspace{0.2cm}

\begin{x}{\small\bf EXAMPLE} \ %08
Let $\rho \in \R_{\alg}$, $\bx = (x_1, \ldots, x_n) \in \R_{\alg}^n$ $-$then
\[
\begin{cases}
\ds\big\{(y_1, \ldots, y_n) \in \R^n \text{:} \ \sum\limits_{i=1}^n \hsx (y_i - x_i)^2 \ < \ \rho^2 \big\}\\[11pt]
\ds\big\{(y_1, \ldots, y_n) \in \R^n \text{:} \ \sum\limits_{i=1}^n \hsx (y_i - x_i)^2 \ \leq \ \rho^2 \big\}
\end{cases}
\]
are semialgebraic as are
\[
\begin{cases}
\ds\big\{(y_1, \ldots, y_n) \in \R^n \text{:} \ \max\limits_{i=1,\ldots, n} \hsx \abs{y_i - x_i}\ < \ \rho \big\}\\[11pt]
\ds\big\{(y_1, \ldots, y_n) \in \R^n \text{:} \ \max\limits_{i=1,\ldots, n} \hsx \abs{y_i - x_i} \ \leq \ \rho \big\}
\end{cases}
.
\]
\end{x}
\vspace{0.3cm}

\begin{x}{\small\bf EXAMPLE} \ %09
\[
\big\{(x,y) \in \R^2 \text{:} \ \exists \ n \in \N, \ y = nx\big\}
\]
is not semialgebraic.
\end{x}
\vspace{0.3cm}

\begin{x}{\small\bf EXAMPLE} \ %10
\[
\big\{(x,y) \in \R^2 \text{:}  \ y = e^x\big\}
\]
is not semialgebraic.
\end{x}
\vspace{0.2cm}

%%----------------------------------------------------------------------------------------------03
\begin{x}{\small\bf DEFINITION} \ %11
\\

\hspace{.9cm}\textbullet \quad A \un{basic open} semialgebraic subset of $\R^n$ is a set of the form
\[
\big\{\bx \in \R^n : \  f_1(\bx) > 0, \ldots, f_r(x) > 0\big\},
\]
where
\[
f_1, \ldots, f_r \in \ \R_{\alg} [T_1, \ldots, T_n].
\]
\\[-.75cm]

\hspace{.9cm}\textbullet \quad A \un{basic closed} semialgebraic subset of $\R^n$ is a set of the form
\[
\big\{\bx \in \R^n : \  f_1(\bx) \geq 0, \ldots, f_r(\bx) \geq 0\big\},
\]
where
\[
f_1, \ldots, f_r \in \ \R_{\alg} [T_1, \ldots, T_n].
\]
\end{x}
\vspace{0.3cm}

\begin{x}{\small\bf LEMMA} \ %12
Suppose that \mX is an open (closed) semialgebraic set $-$then \mX is a finite union of basic open (basic closed) semialgebraic sets.
\end{x}
\vspace{0.3cm}

\begin{x}{\small\bf DEFINTION} \ %13
\ Let $X \subset \hsx \R^n$, $Y \subset \hsx \R^m$ be semialgebraic sets $-$then a \un{semialgebraic map}
$f:X \ra Y$ is a continuous function such that the graph 
$\Gamma_f \subset \ X \times Y$ is a semialgebraic subset of $\R^n \times \R^m$.
\end{x}
\vspace{0.3cm}

\begin{x}{\small\bf SCHOLIUM} \ %14
Let
\[
P_1, \ldots, P_m
\]
be elements of $\R_{\alg}[T_1, \ldots, T_n]$ $-$then the arrow $f:X \ra Y$ defined by the prescription 
\[
(x_1, \ldots, x_n) \mapsto (P_1(x_1, \ldots, x_n), \ldots, P_m(x_1, \ldots, x_n))
\]
is a semialgebraic map.
\\

[Note: \ 
One can replace $P_1, \ldots, P_m$ by elements 
\[
\frac{P_1}{ Q_1}, \ldots, \frac{P_m}{ Q_m}
\]
%%----------------------------------------------------------------------------------------------04
 of $\R_{\alg}(T_1, \ldots, T_m)$ provided that none of the 
 \[
 Q_1, \ldots, Q_m
 \]
 vanish at any point of \mX.]
\end{x}
\vspace{0.3cm}

\begin{x}{\small\bf EXAMPLE} \ %15
The graph of $f(x) = \sqrt{x}$ is the semialgebraic set 
\[
\big\{(x,y) \in \ \R^2 \text{:} \ y^2 = x, y \geq 0\big\},
\]
so $f$ is a semialgebraic map.
\end{x}
\vspace{0.3cm}

\begin{x}{\small\bf EXAMPLE} \ %16
If \mX is semialgebraic, then the diagonal $\Delta:X \ra X \times X$ is a semialgebraic map.

[The graph of $\Delta$ is the intersection of $X \times X \times X$ with the semialgebraic set
\[
\big\{(x,y,z) \text{:} \ x = y = z \big\}.]
\]
\end{x}
\vspace{0.3cm}

\begin{x}{\small\bf EXAMPLE} \ %17
If $\emptyset \neq X \subset \ \R^n$ is a semialgebraic set, then the distance function 
$x \ra \text{dist}(x,X)$ is a semialgebraic map.
\end{x}
\vspace{0.3cm}

\begin{x}{\small\bf DEFINTION} \ %18
Let $X \subset \hsx \R^n$ be an open semialgebraic subset $-$then an analytic function $\phi : X \ra \R$ is said to satisfy the 
\un{condition of Nash} if there are elements 
\[
a_0, a_1, \ldots, a_d \ \neq \ 0
\]
of $\R_{\alg}[T_1, \ldots, T_n]$ such that
\[
a_0 + a_1 \phi + \cdots + a_d \phi^d \ = \ 0,
\]
i.e., $\phi$ is an analytic algebraic function.  
Such a function is necessarily a semialgebraic map.
\end{x}
\vspace{0.3cm}

\begin{x}{\small\bf EXAMPLE} \ %19
Let $X = ]-1,1[ \ \subset \ \R$ and take 
$\phi(x) = \sqrt{1 - x^2}$ $-$then
\[
\phi(x)^2 \ - \ (1 - x^2) \ = \ 0,
\]
%%----------------------------------------------------------------------------------------------05
so $\phi$ satisfies the condition of Nash.
\end{x}
\vspace{0.3cm}

\begin{x}{\small\bf LEMMA} \ %20
Let $X \subset \ \R^n$, $X^\prime \subset \ \R^{n^\prime}$, $Y \subset \ \R^m$, $Y^\prime \subset \ \R^{m^\prime}$ 
be semialgebraic sets and let $f:X \ra Y$, $f^\prime : X^\prime \ra Y^\prime$ be semialgebraic maps $-$then
\[
f \times f^\prime : X \times X^\prime \ra Y \times Y^\prime
\]
is a semialgebraic map.
\\

PROOF \ 
The graph of $f \times f^\prime$ is
\[
\big\{(x,x^\prime, y, y^\prime) \ \text{:} \ y = f(x), \ y^\prime = f^\prime(x^\prime)\big\}
\]
which is the intersection of 
\[
\big\{(x,x^\prime, y, y^\prime) \ \text{:} \ y = f(x)\big\} \approx \Gamma_f \times \R^{n^\prime} \times \R^{m^\prime}
\]
and 
\[
\big\{(x,x^\prime, y, y^\prime) \ \text{:} \ y^\prime = f^\prime(x^\prime)\big\} \approx \Gamma_{f^\prime} \times \R^{n} \times \R^{m}.
\]
\end{x}
\vspace{0.3cm}

\begin{x}{\small\bf THEOREM} \ (Tarski-Seidenberg) \ %21
If $f : \R^n \ra \R^m$ is a polynomial function (cf. \#14) and $X \subset \ \R^n$ is a semialgebraic set, then $f(X)$ is semialgebraic.
\end{x}

\begin{x}{\small\bf SCHOLIUM} \ %22
Let $X \subset \ \R^{n+m}$ be semialgebraic, $\Pi : \R^{n+m} \ra \R^n$ the projection onto the space of the first $n$ coordinates 
(or $\Pi : \R^{n+m} \ra \R^m$ the projection onto the space of the second $m$ coordinates) $-$then $\Pi(X)$ is a semialgebraic subset of 
$\R^n$ (or $\R^m$).
\end{x}
\vspace{0.3cm}

\begin{x}{\small\bf LEMMA} \ %23
Let $X \subset \ \R^n$, $Y \subset \ \R^m$ be semialgebraic sets and let $f:X \ra Y$ be a semialgebraic map $-$then $f(X)$ is a semialgebraic set.
\\

PROOF \ 
$\Gamma_f$ is a semialgebraic set by assumption.  Now apply \#22 to
\[
\Gamma_f \subset \ \R^n \ \times \R^m \ \equiv \ \R^{n+m} \ \ra \ \R^m.
\]
\end{x}
\vspace{0.3cm}
%%----------------------------------------------------------------------------------------------06

\begin{x}{\small\bf APPLICATION} \ %24
If $X \subset \ \R^n$, $Y \subset \ \R^m$, $Z \subset \ \R^{\ell}$ are semialgebraic sets and if 
$f:X \ra Y$, $g:Y \ra Z$  are semialgebraic maps, then $g \circ f$ is a semialgebraic map.

[In fact,
\[
\Gamma_{g \circ f} 
\ = \ 
(1_X \times g) \hsx (\Gamma_f).]
\]
\end{x}
\vspace{0.3cm}

\begin{x}{\small\bf \un{N.B.}} \ %25
It follows that there is a category whose objects are the semialgebraic sets and whose morphisms are the semialgebraic maps.
\end{x}

\begin{x}{\small\bf REMARK} \ %26
Let \mS, \mX, \mY be semialgebraic sets and let $f:X \ra S$, $g:Y \ra S$ be semialgebraic maps $-$then the fiber product 
$X \times_S Y$ of sets is a semialgebraic set and is a fiber product in the category of semialgebraic sets and semialgebraic maps.
\end{x}
\vspace{0.3cm}

\begin{x}{\small\bf LEMMA} \ %27
The inverse image of a semialgebraic set under a semialgebraic map is a semialgebraic set.
\\[-.25cm]

PROOF \ 
Let $X \subset \ \R^n$, $Y \subset \ \R^m$, $Z \subset \ Y$ be semialgebraic sets and suppose that $f:X \ra Y$ is a semialgebraic map.  
Write
\[
f^{-1}(Z) \ = \ \Pi(\Gamma_f \ \cap \ (\R^n \times Z))
\]
and apply \#22.
\end{x}
\vspace{0.3cm}

\begin{x}{\small\bf \un{N.B.}} \ %28
The converse of this lemma is false.

[The exponential function $\exp : \R \ra \R$ is not semialgebraic but does have the property that the inverse image of a semialgebraic set is a semialgebraic set.]
\end{x}

\begin{x}{\small\bf THEOREM} \ %29
Suppose that \mX is a semialgebraic set $-$then the semialgebraic maps $X \ra \R$ form a ring under pointwise addition and multiplication.
\end{x}
\vspace{0.3cm}

%%----------------------------------------------------------------------------------------------07

\begin{x}{\small\bf LEMMA} \ %30
If $X \subset \hsx \R^n$ is a semialgebraic set, then the closure, the interior, and the frontier of \mX are semialgebraic.
\end{x}
\vspace{0.3cm}

\begin{x}{\small\bf EXAMPLE} \ %31
Take
\[
X \ = \ 
\big\{x \in \R \text{:} \ (x^2 - 1) (x - 2)^2 < 0\big\}.
\]
Then $\ov{X} = [-1,1]$ but 
\[
\big\{x \in \R \text{:} \ (x^2 - 1) (x - 2)^2 \leq 0\big\} 
\ = \ 
\ov{X} \cup \{2\}.
\]
\end{x}
\vspace{0.3cm}

%%%%%%%%%%%%%%%%%%%%%%%%%%%%%%%%%%%%%%
%%%%%%%%%%%%%%%%%%%%%%%%%%%%%%%%%%%%%%
%%%%%%%%%%%%%%%%%%%%%%%%%%%%%%%%%%%%%%

\begingroup
%\fontsize{9pt}{11pt}\selectfont

\begin{center}
%$\S \ 0$\\

$\mathcal{REFERENCES}$\\
\vspace{0.5cm}
%\textbf{BOOKS}\\
\end{center}

\setlength\parindent{0 cm}

[1] \quad Benedetti, R., Risler, J-J., \textit{Real Algebraic and Semialgebraic Sets}, Hermann (1990).\\

[2] \quad Bochnak, J., Coste, M., Roy, M-F., \textit{Real Algebraic Geometry}, Springer (1998).\\

%123
\chapter{
$\boldsymbol{\S}$\textbf{2}.\quad  PERIODS: DEFINITION AND EXAMPLES}
\setlength\parindent{2em}
\setcounter{theoremn}{0}
\renewcommand{\thepage}{\S2-\arabic{page}}
%%----------------------------------------------------------------------------------------------01

\begin{x}{\small\bf NOTATION} \ %01
$\sS A^n$ is the set of semialgebraic subsets of $\R^n$ with a nonempty interior.
\end{x}

\begin{x}{\small\bf DEFINITION} \ %02
A real number $\pe$ is a \un{period} in the sense of Kontsevich-Zagier if

\hspace{0.5cm} \textbullet \quad $\exists \ n \in \ \N$
\\[-.25cm]

\hspace{0.5cm} \textbullet \quad $\exists \ X \in \ \sS A^n$
\\[-.25cm]

\hspace{0.5cm} \textbullet \quad $\exists \ P, \ Q \in \ \R_{\alg}[T_1, \ldots, T_n] \ (\restr{Q}{X} \neq 0) $
\\[.25cm]
such that 
\begin{align*}
\pe \ 
&= \ 
\tI\bigg(X,\frac{P}{Q}\bigg) 
\\[11pt]
&\equiv \ 
\ds \int\limits_X \hsx \frac{P(x_1, \ldots, x_n)}{Q(x_1, \ldots, x_n)}  \hsx \td x_1 \cdots \td x_n
\\[11pt]
&\equiv \ 
\int\limits_{X} \hsx \frac{P}{Q} (\bx) \td \bx
\end{align*}
is an absolutely convergent integral, $\td \bx = \td x_1 \cdots \td x_n$ being Lebesgue measure.
\end{x}

\begin{x}{\small\bf NOTATION} \ %03
Write
\[
\sP_\text{KZ}
\]
for the set of periods in the sense of Kontsevich-Zagier.

[Note: \ 
We shall work exclusively in the real domain, a period in the complex domain being a combination 
$\pe_1 + \sqrt{-1} \hsx \pe_2$, where $\pe_1, \ \pe_2 \in \sP_\text{KZ}$.]
\end{x}
\vspace{0.2cm}

\begin{x}{\small\bf LEMMA} \ %04
\[
\R_\alg \ \subset \ \sP_\text{KZ}.
\]
%%----------------------------------------------------------------------------------------------02

[Take $\rho \in \ \R_\alg$, $\rho > 0$, $-$then 
\[
\rho 
\ = \ 
\int\limits_{0 < x < \rho} \hsx 1 \hsx \td x.]
\]
\end{x}
\vspace{0.2cm}

\begin{x}{\small\bf EXAMPLE} \ %05
Suppose that  $\rho \in \ \R_\alg$, $\rho > 1$, $-$then 
\[
\elln(\rho) 
\ = \ 
\int_1^\rho \hsx \frac{\td t}{t} \ \in \ \sP_\text{KZ}.
\]

[Note: \ 
Observe too that 
\[
\int_1^\rho  \hsx \frac{\td t}{t} 
\ = \ 
\int\limits_{\substack{1 < x < \rho\\[2pt] 0 < x y < 1}} \hsx 1 \hsx \td x \td y.]
\]
\end{x}
\vspace{0.3cm}

\begin{x}{\small\bf EXAMPLE} \ %06
Consider
\[
\int\limits_{x^2 + y^2 \leq 1} \hsx 1 \hsx \td x \td y
\]
to see that $\pi \in \ \sP_\text{KZ}$.
\\[-.25cm]

[It is unknown whether $\ds \frac{1}{\pi} \in \ \sP_\text{KZ}$ or not.]
\end{x}
\vspace{0.3cm}

\begin{x}{\small\bf EXAMPLE} \ %07
\begin{align*}
\zeta(2) \ 
&=\ 
\int\limits_{1 > x_1 > x_2 > 0} \hsx \frac{\td x_1}{x_1} \cdot \frac{\td x_2}{1 - x_2}
\\[11pt]
&\in \ 
\sP_\text{KZ}.
\end{align*}
\end{x}
\vspace{0.2cm}

\begin{x}{\small\bf EXAMPLE} \ %08
Take for \mS the square
\[
\begin{cases}
\ 0 \ \leq x \ \leq \ 1\\
\ 0 \ \leq y \ \leq \ 1
\end{cases}
.
\]
Let 
\[
f(x,y) 
\ = \ 
\frac{x y}{(x^2 + y^2)^2}
\]
%%----------------------------------------------------------------------------------------------03
if $x^2 + y^2 > 0$ and set $f(0,0) = 0$ $-$then $f$ is not integrable.  
For if it were, then the iterated integral
\[
\int_0^1 \hsx \td x \ \int_0^1 \frac{xy}{(x^2 + y^2)^2} \hsx \td y
\]
would exist (Fubini).  
But for $x \neq 0$, 
\[
\int_0^1 \frac{xy}{(x^2 + y^2)^2} \hsx \td y 
\ = \ 
\frac{1}{2 x} \ - \ \frac{x}{2(x^2 + 1)}
\]
and this function is not integrable in $]0,1]$.
\\[-.25cm]

[Note: \ 
Replace
\[
\frac{xy}{(x^2 + y^2)^2} \qquad \text{by} \qquad \frac{x^2 - y^2}{(x^2 + y^2)^2} \hsx.
\]
Then
\[
\begin{cases}
\ \ds\int_0^1 \hsx \td x \ \int_0^1 \hsx f(x,y) \td y \ = \ \ \ \frac{\pi}{4} \\[15pt]
\ \ds \int_0^1 \hsx \td y \ \int_0^1 \hsx f(x,y) \td x \ = \ -\frac{\pi}{4}
\end{cases}
.]
\]
\end{x}
\vspace{0.3cm}

%%----------------------------------------------------------------------------------------------07
%%%%%%%%%%%%%%%%%%%%%%%%%%%%%%%%%%%%%%
%%%%%%%%%%%%%%%%%%%%%%%%%%%%%%%%%%%%%%
%%%%%%%%%%%%%%%%%%%%%%%%%%%%%%%%%%%%%%

%123
\chapter{
$\boldsymbol{\S}$\textbf{3}.\quad  PERIODS: ALGEBRAIC CONSIDERATIONS}
\setlength\parindent{2em}
\setcounter{theoremn}{0}
\renewcommand{\thepage}{\S3-\arabic{page}}
%%----------------------------------------------------------------------------------------------01

\ \indent 
It was pointed out in $\S 2$, \#4 that
\[
\R_\alg \ \subset \ \sP_\text{KZ}.
\]
Of course $\R_\alg$ is a countable field.  And:

\begin{x}{\small\bf LEMMA} \ %01
$\sP_\text{KZ}$ is countable.
\\

PROOF \ 
$\forall \ n \in \N$,  $\R_\alg[T_1, \ldots, T_n]$ is countable, hence $\R_\alg(T_1, \ldots, T_n)$ is countable,  
as is $\sS A^n$, hence $\sP_\text{KZ}$ is countable.

[Note: \ 
Consequently ``most'' real numbers are not periods.  
And if a real number is not a period, then it is transcendental.]
\end{x}

\begin{x}{\small\bf REMARK} \ %02
It is unknown if $\sP_\text{KZ}$ is a field (but $\sP_\text{KZ}$ is an $\R_\alg$ algebra (see \#8 below)).
\end{x}

\begin{x}{\small\bf DISCUSSION} \ %03
Suppose that
\[
\pe 
\ = \ 
I\bigg(X,\frac{P}{Q}\bigg) 
\ = \ 
\int_X \hsx \frac{P}{Q}(\bx) \td \bx \qquad \text{(cf. $\S 2$, $\#2$)}
\]
and
\[
\int_X \hsx \abs{\frac{P}{Q}(\bx)} \td \bx \ < \ \infty,
\]
hence 
\[
\int_X \hsx \abs{\frac{P}{Q}(\bx)} \td \bx \ = \ 0
\]
iff
\[
\abs{\frac{P}{Q}(\bx)} \ = \ 0
\]
%%----------------------------------------------------------------------------------------------02
almost everywhere, thereby forcing $\pe = 0$.
\end{x}
\vspace{0.2cm}

\begin{x}{\small\bf NOTATION} \ %04
Given a measurable subset $X \subset \ \R^n$, let
\begin{align*}
\vol_n(X) \ 
&= \ 
\int_X \hsx 1 \hsx \td x_1 \ldots \hsx \td x_n
\\[15pt]
&\equiv \ 
\int_X \hsx 1 \hsx \td \bx.
\end{align*}
\end{x}
\vspace{0.2cm}

\begin{x}{\small\bf EXAMPLE} \ %05
Consider the $n$-simplex $\Delta_n$ $-$then
\[
\vol_n(\Delta_n) 
\ = \ 
\frac{1}{n!}.
\]
\end{x}
\vspace{0.3cm}

\begin{x}{\small\bf EXAMPLE} \ %06
Consider the $n$-ball $\sB^n$ $-$then
\[
\vol_n(\sB^n) 
\ = \ 
\frac{\pi^{\frac{n}{2}}}{\Gamma\bigg(\frac{n}{2} + 1\bigg)} \hsx .
\]
\end{x}
\vspace{0.3cm}

\begin{x}{\small\bf LEMMA} \ %07
Let $\pe \in \R$ $-$then $\pe \in \sP_\text{KZ}$ iff for some $d \in \ \N$, there exist disjoint semialgebraic sets 
$X_1 \subset \ \R^d$, $X_2 \subset \ \R^d$ of finite volume such that
\[
\pe \ = \ \vol_d(X_1) - \vol_d(X_2).
\]

PROOF \ 
It can be assumed that $\pe \neq 0$.

\hspace{0.5cm} \textbullet \quad Given $\pe \in \sP_\text{KZ}$, put
\[
\begin{cases}
\ \ds X_+ \ = \ \big\{\bx \in \ X \text{:} \ \sgn \frac{P}{Q} (\bx) \ = \ +1\big\}\\[15pt]
\ \ds X_- \ = \ \big\{\bx \in \ X \text{:} \ \sgn \frac{P}{Q} (\bx) \ = \ -1\big\}
\end{cases}
,
\]
a disjoint partition of \mX (to within a set of measure 0) (cf. \#3) $-$then
\[
\pe
\ = \ 
I\bigg(X,\frac{P}{Q}\bigg) 
\ = \ 
I\bigg(X_+,\frac{P}{Q}\bigg) - I\bigg(X_-,-\frac{P}{Q}\bigg).
\]
%%----------------------------------------------------------------------------------------------03
Now introduce semialgebraic sets
\[
\begin{cases}
\ \ds X_1 \ = \ \big\{(\bx,t) \in \ X \times \R \text{:} \ t > 0, t \leq \frac{P}{Q} (\bx)\big\}\\[15pt]
\ \ds X_2 \ = \ \big\{(\bx,t) \in \ X \times \R \text{:} \ t < 0, t \geq \frac{P}{Q} (\bx)\big\}
\end{cases}
\]
from which
\[
\begin{cases}
\ \ds I\bigg(X_+,\frac{P}{Q}\bigg) \ = \ \int_{X_1} \hsx 1 \hsx \td \bx \td t \\[15pt]
\ \ds I\bigg(X_-,-\frac{P}{Q}\bigg) \ = \ \int_{X_2} \hsx 1 \hsx \td \bx \td t 
\end{cases}
.
\]
Therefore
\begin{align*}
\pe \ 
&= \ 
I\bigg(X,\frac{P}{Q}\bigg) 
\\[15pt]
&= \ 
\int_{X_1} \hsx 1 \hsx \td \bx \td t  -  \int_{X_2} \hsx 1 \hsx \td \bx \td t 
\\[15pt]
&= \ 
\vol_{n+1} (X_1) - \vol_{n+1} (X_2).
\end{align*}
Matters are thus settled with the choice $d = n + 1$.
\\[-.25cm]

\hspace{0.5cm} \textbullet \quad Suppose that
\[
\pe \ = \ \vol_d(X_1) - \vol_d(X_2),
\]
the claim being that $\pe \in \  \sP_\text{KZ}$.  
To see this, write
\begin{align*}
\pe \ 
&= \ 
\int_{X_1} \hsx 1 \td x_1 \cdots \td x_d - \int_{X_2} \hsx 1 \td x_1 \cdots \td x_d
\\[15pt]
&= \ 
\int_{X_1} \hsx \bigg(\int_0^1 \hsx 2 t \td t\bigg) \td x_1 \cdots \td x_d 
+ 
\int_{X_2} \hsx \bigg(\int_{-1}^0 \hsx 2 t \td t\bigg) \td x_1 \cdots \td x_d 
\\[15pt]
&= \ 
\int_{Y_1 \cup Y_2} \hsx \hsx 2 t \hsx \td t\td x_1 \cdots \td x_d,
\end{align*}
where
\[
Y_1 \ = \ ]0,1[ \ \times \ X_1, \quad
Y_2 \ = \ ]-1,0[ \ \times X_2
\]
%%----------------------------------------------------------------------------------------------04
are disjoint semialgebraic sets in $\R^1 \times \R^d$.  
Therefore $\pe \in \sP_\text{KZ}$.
\end{x}
\vspace{0.2cm}

\begin{x}{\small\bf THEOREM} \ %08
$\sP_\text{KZ}$ is an $\R_\alg$ algebra.
\\

PROOF \ 
There are two issues
\[
\begin{cases}
\ \text{Stability of the product of two periods}\\[8pt]
\ \text{Stability of the sum of two periods}
\end{cases}
.
\]

\hspace{0.5cm} \textbullet \quad (Product) \ 
Given
\[
\pe_1 \ = \ \int_{X_1} \hsx \frac{P_1}{Q_1} (\bx_1) \hsx \td \bx_1, 
\qquad
\pe_2 \ = \ \int_{X_2} \hsx \frac{P_2}{Q_2} (\bx_2) \hsx \td \bx_2,
\]
write
\[
\begin{cases}
\ \pro_1 \text{:} \R^{d_1} \times \R^{d_2} \ra \R^{d_1} \\[8pt]
\ \pro_2 \text{:} \R^{d_1} \times \R^{d_2} \ra \R^{d_2}
\end{cases}
\]
and define
\[
F \ \text{:} \ \R^{d_1} \times \R^{d_2}  \ra \R
\]
by
\[
F(\bx_1, \bx_2) 
\ = \ 
\frac{P_1}{Q_1} \circ \pro_1(\bx_1, \bx_2) 
\cdot 
\frac{P_2}{Q_2} \circ \pro_2(\bx_1, \bx_2).
\]
Then
\begin{align*}
\int\limits_{X_1 \times X_2} \hsx F(\bx_1, \bx_2) \hsx \td \bx_1 \td \bx_2 \ 
&=\ 
\bigg(\int\limits_{X_1} \frac{P_1}{Q_1} (\bx_1) \hsx \td \bx_1 \bigg) \hsx
\bigg(\int\limits_{X_2} \frac{P_2}{Q_2} (\bx_2) \hsx \td \bx_2 \bigg) 
\\[15pt]
&=\ 
\pe_1 \pe_2.
\end{align*}
%%----------------------------------------------------------------------------------------------05
\hspace{0.5cm} \textbullet \quad (Sum) \quad
Let $\pe_1$, $\pe_2 \in \  \sP_\text{KZ}$.  
Per $\# 7$, write 
 \[
\begin{cases}
\ \pe_1 \ = \  \vol_{d_1} (X_1) - \vol_{d_1} (X_2) \\[8pt]
\ \pe_2 \ = \  \vol_{d_2} (Y_1) - \vol_{d_2} (Y_2) 
\end{cases}
,
\]
where
\[
\begin{cases}
\ X_1 \ \subset \ \R^{d_1}, \ X_2 \ \subset \ \R^{d_1} \\[8pt]
\ Y_1 \ \subset \ \R^{d_2}, \ Y_2 \ \subset \ \R^{d_2}
\end{cases}
\]
are semialgebraic and 
\[
X_1 \ \cap \ X_n \ = \ \emptyset, \qquad 
Y_1 \ \cap \ Y_2 = \emptyset.
\]
There is no loss of generality in assuming that $d_1 = d_2 = d$ 
(if, e.g., $d_1 < d_2$, let $k = d_2 - d_1$ and work with 
$X_1 \times [0,1]^k, \hsx X_2 \times [0,1]^k$).  
This said, it then follows that
\[
\pe_1 + \pe_2 
\ = \ 
\vol_d(X_1) + \vol_d(Y_1) - \vol_d(X_2) - \vol_d(Y_2)
\]
or still, 
\[
\pe_1 + \pe_2 
\ = \ 
\vol_d(X_1 \cup Y_1) + \vol_d(X_1 \cap Y_1) - \vol_d(X_2 \cup Y_2) - \vol_d(X_2 \cap Y_2).
\]
Put
\[
\begin{cases}
\ 
\ W_1 \ = \ X_1 \cup Y_1, \quad W_2 \ = \ X_2 \cup Y_2 \\[8pt]
\ Z_1 \ \ \ = \ X_1 \cap Y_1, \quad Z_2 \ \ = \ X_2 \cap Y_2
\end{cases}
\]
and let $I_1, J_1, I_2, J_2 \subset \ \R$ be disjoint open intervals of length 1 with endpoints in $\R_\alg$ $-$then
\[
\pe_1 + \pe_2 \ = \ \vol_{d+1} ((W_1 \times I_1) \cup Z_1 \times J_1))
\ - \ 
\vol_{d+1}((W_2 \times I_2) \cup (Z_2 \times J_2)).
\]
%%----------------------------------------------------------------------------------------------06
Here
\[
((W_1 \times I_1) \cup (Z_1 \times J_1)) \cap ((W_2 \times I_2) \cup (Z_2 \times J_2)) \ = \ \emptyset.
\]
In addition
\[
\begin{cases}
\ (W_1 \times I_1) \cup (Z_1 \times J_1) \\[8pt]
\ (W_2 \times I_2) \cup (Z_2 \times J_2)
\end{cases}
\]
are semialgebraic subsets of $\R^{d+1}$, thus to finish it remains only to quote $\# 7$.
\end{x}
\vspace{0.3cm}

\begin{x}{\small\bf APPLICATION} \ %09
Let $x \in \ \R_\alg^\times$, $\pe \in \ \sP_\text{KZ}$ $-$then
\[
\begin{cases}
\ x + \pe \in \ \sP_\text{KZ}\\[8pt]
\ \hspace{.5cm} x \pe \in \ \sP_\text{KZ}
\end{cases}
.
\]
\end{x}
\vspace{0.3cm}

%%----------------------------------------------------------------------------------------------07
%%%%%%%%%%%%%%%%%%%%%%%%%%%%%%%%%%%%%%
%%%%%%%%%%%%%%%%%%%%%%%%%%%%%%%%%%%%%%
%%%%%%%%%%%%%%%%%%%%%%%%%%%%%%%%%%%%%%

%123
\chapter{
$\boldsymbol{\S}$\textbf{4}.\quad  PERIODS: THEORETICAL CONSIDERATIONS}
\setlength\parindent{2em}
\setcounter{theoremn}{0}
\renewcommand{\thepage}{\S4-\arabic{page}}
%%----------------------------------------------------------------------------------------------01

\begin{x}{\small\bf THEOREM} \ %01
Let
\[
\pe
\ = \ 
I\bigg(X,\frac{P}{Q}\bigg) \qquad (X \in \ \sS A^n)
\]
be a nonzero period $-$then there are compact semialgebraic sets 
$K_1, \ldots, K_m \ in \ \sS A^n$, polynomials 
\[
\begin{cases}
\ P_1, \ldots, P_m \in \ \R_\alg[T_1, \ldots, T_n]\\[8pt]
\ Q_1, \ldots, Q_m \in \ \R_\alg[T_1, \ldots, T_n]
\end{cases}
\]
with $\restr{Q_1}{K_1} \neq 0, \ldots, \restr{Q_m}{K_m} \neq 0$ such that
\[
I\bigg(X,\frac{P}{Q}\bigg)
\ = \ 
\sum\limits_{i=1}^m \ I\bigg(K_i,\frac{P_i}{Q_i}\bigg).
\]

[This result is due to Juan Viu-Sos 
\footnote[2]{https://arxiv.org/abs/1509.01097}. %\footnote[2]{\vspace{.11 cm} arXiv: \ 1509.01097 [math,NT].}
Its proof is difficult, depending, as it does, on Hironaka's rectilinearization of semialgebraic sets.]
\end{x}

\begin{x}{\small\bf \un{N.B.}} \ %02
The integrals
\[
I\bigg(K_i,\frac{P_i}{Q_i}\bigg)
\]
are absolutely convergent.
\end{x}

\begin{x}{\small\bf LEMMA} \ %03
Let
\[
\pe 
\ = \ 
I\bigg(X,\frac{P}{Q}\bigg) \qquad (X \in \ \sS A^n)
\]
be a nonzero period $-$then there are compact semialgebraic sets \mS and \mT in $\sS A^{n+1}$
%%----------------------------------------------------------------------------------------------02
such that 
\[
\pe 
\ = \ 
\vol_{n+1}(S) - \vol_{n+1}(T).
\]

PROOF \ 
Proceed as in $\S 3$, $\# 7$ (necessity), thus 
\[
I\bigg(X,\frac{P}{Q}\bigg)
\ = \ 
I\bigg(X_+,\frac{P}{Q}\bigg) - I\bigg(X_-,-\frac{P}{Q}\bigg).
\]
Per $\# 1$, write  %%%dmcXXX do we need minus signs near the end of both the next lines
\[
\begin{cases}
\ \ds I\bigg(X_+,\frac{P}{Q}\bigg) \hspace{0.4cm} = \ \sum\limits_{i=1}^{m_+} \ I\bigg(K_i^+,\frac{P_i^+}{Q_i^+}\bigg)\\[21pt]
\ \ds I\bigg(X_-,-\frac{P}{Q}\bigg) \ = \ \sum\limits_{i=1}^{m_-} \ I\bigg(K_i^-,\frac{P_i^-}{Q_i^+-}\bigg)
\end{cases}
\]
or still, 
\[
\begin{cases}
\ \ds I\bigg(X_+,\frac{P}{Q}\bigg) \hspace{0.4cm}  = \ \sum\limits_{i=1}^{m_+} \ \int\limits_{L_i^+} \hsx 1 \hsx \td \bx \td t\\[21pt]
\ \ds I\bigg(X_-,-\frac{P}{Q}\bigg) \ = \ \sum\limits_{i=1}^{m_-} \ \int\limits_{L_i^-} \hsx 1 \hsx \td \bx \td t
\end{cases}
\]
or still, 
\[
\begin{cases}
\ \ds I\bigg(X_+,\frac{P}{Q}\bigg) \hspace{0.4cm}  = \ \sum\limits_{i=1}^{m_+} \ \vol_{n+1} (L_i^+)\\[21pt]
\ \ds I\bigg(X_-,-\frac{P}{Q}\bigg) \ = \ \sum\limits_{i=1}^{m_-} \ \vol_{n+1} (L_i^-)
\end{cases}
,
\]
where $L_i^+$, $L_i^- \in \ \sS A^{n+1}$ are compact.  
Working with semialgebraic translations, it can be arranged that the $L_i^+$ are pairwise disjoint 
and the $L_i^-$ are pairwise disjoint.
%%----------------------------------------------------------------------------------------------03
Put now
\[
S 
\ = \ 
\bigcup\limits_{i=1}^{m_+} \ L_i^+, 
\quad
T 
\ = \ 
\bigcup\limits_{i=1}^{m_-} \ L_i^-.
\]
Then
\begin{align*}
\pe \ 
&= \ 
I\bigg(X,\frac{P}{Q}\bigg)
\\[15pt]
&= \
\sum\limits_{i=1}^{m_+} \ \vol_{n+1} (L_i^+) - \sum\limits_{i=1}^{m_-} \ \vol_{n+1} (L_i^-)
\\[15pt]
&= \
\vol_{n+1} \bigg(\bigcup\limits_{i=1}^{m_+} \ L_i^+\bigg) - \vol_{n+1} \bigg(\bigcup\limits_{i=1}^{m_-} \ L_i^-\bigg)
\\[15pt]
&= \
\vol_{n+1} (S) - \vol_{n+1} (T).
\end{align*}
\end{x}
\vspace{0.2cm}

\begin{x}{\small\bf THEOREM} \ %04
Let $\pe \in \sP_\text{KZ}$ $(\pe \neq 0)$, say
\[
\pe 
\ = \ 
I\bigg(X,\frac{P}{Q}\bigg) \qquad (X \in \ \sS A^n).
\]
Then there exists a compact $K \in \ \sS A^k$ $(0 < k \leq n + 1)$ such that 
\[
\pe \ = \ \sgn(\pe) \cdot \vol_k(K).
\]

It can be assumed that $\pe$ is positive and in the notation of $\# 3$, 
\[
0 < \vol_{n+1} (T) \ < \ \vol_{n+1}(S).
\]
The point then is to construct a compact semialgebraic set \mK from \mS and \mT so as to arrive at
\[
\pe \ = \ \vol_k(K) \qquad (0 < k \leq n+1).
\]
While ``elementary'', the details are tedious and will be omitted.
\end{x}
\vspace{0.2cm}
%%----------------------------------------------------------------------------------------------04

\begin{x}{\small\bf EXAMPLE} \ %05
$\pi^2$ is the 4-dimensional volume of the product of two copies of the unit disk and the 3-dimensional volume of the set
\[
\big\{(x,y,z) \in \ \R^3 \ \text{:} \ x^2 + y^2 \leq 1, \ 0 \leq z \ ((x^2 + y^2)^2 + 1) \leq 4\big\}.
\]
\end{x}
\vspace{0.3cm}

\begin{x}{\small\bf DEFINITION} \ %06
Let $\pe \in \sP_\text{KZ}$ $(\pe \neq 0)$ $-$then the \un{degree} of $\pe$, denoted $\deg(\pe)$, 
is the smallest positive integer $k$ such that
\[
\abs{\pe} 
\ = \ 
\vol_k(K) \qquad (K \in \ \sS A^k, K \ \text{compact) (cf. $\#4$)}.
\]

[Note: \ 
Take $\deg(0) = 0$ and in the complex domain, let
\[
\deg(\pe_1 + \sqrt{-1} \hsx \pe_2) 
\ = \ 
\max \{\deg(\pe_1), \deg(\pe_2)\}.]
\]
\end{x}
\vspace{0.3cm}

\begin{x}{\small\bf EXAMPLE} \ %07
\[
\deg(\pi) = 2.
\]
[In view of $\S 2$, $\#6$, $\deg(\pi) \leq 2$.  
On the other hand, the fact that $\pi$ is transcendental rules out the possibility that
\[
\deg(\pi)  \ = \ 1 \qquad \text{(cf. $\# 9$ infra).}]
\]
\end{x}
\vspace{0.2cm}

\begin{x}{\small\bf REMARK} \ %08
It is conjectured that $\forall \ n \in \ \N$, 
\[
\deg(\pi^n) \ = \ n + 1.
\]
\end{x}
\vspace{0.1cm}

\begin{x}{\small\bf LEMMA} \ %09
$\deg(\pe) = 1$ iff $\pe \in \ \R_\alg^\times$.
\\[-.25cm]

[Since \mK is compact semialgebraic, it can be written as a finite disjoint union of points and open intervals.  
In the other direction, any nonzero $\rho \in \ \R_\alg$ can be written up to sign as the length of $[0,\rho]$ 
(cf. $\S 2$, $\#4$.]
\end{x}
\vspace{0.3cm}
%%----------------------------------------------------------------------------------------------05

\begin{x}{\small\bf APPLICATION} \ %10
A period $\pe \in \sP_\text{KZ}$ is transcendental iff $\deg(\pe) \geq 2$.
\end{x}
\vspace{0.1cm}

\begin{x}{\small\bf LEMMA} \ %11
Let $\pe_1, \ \pe_2 \in  \sP_\text{KZ}$ $-$then
\\

\hspace{0.5cm}  \textbullet \quad 
$\deg(\pe_1 + \pe_2) \leq \max\{\deg(\pe_1), \deg(\pe_2)\}$

\noindent and

\hspace{0.5cm}  \textbullet \quad 
$\deg(\pe_1 \pe_2) \leq \deg(\pe_1) + \deg(\pe_2)$.

\end{x}
\vspace{0.3cm}

\begin{x}{\small\bf RAPPEL} \ %12
$\sP_\text{KZ}$ is an $\R_\alg$ algebra (cf. $\S 3$, $\#8$).
\end{x}
\vspace{0.3cm}

\begin{x}{\small\bf LEMMA} \ %13
Let $x \in \R_\alg^\times$, $\pe \in \sP_\text{KZ}$.  
Assume $\pe \not\in \R_\alg$ $-$then
\[
\deg(x + \pe) \ = \ \deg(\pe).
\]

[In fact, 
\begin{align*}
\deg(\pe) \ 
&=\ 
\deg(-x + x + \pe)
\\[15pt]
&\leq\ 
\max\{\deg(-x), \deg(x + \pe)\}
\\[15pt]
&=\ 
\max\{1, \deg(x+ \pe)\}
\\[15pt]
&=\ 
\deg(x+\pe)
\\[15pt]
&\leq\ 
\max\{1, \deg(\pe)\}
\\[15pt]
&=\ 
\deg(\pe).]
\end{align*}
\end{x}
\vspace{0.3cm}

\begin{x}{\small\bf LEMMA} \ %14
Let $x \in \ \R_\alg^\times$, $\pe \in \sP_\text{KZ}$.  
Assume: $\pe \not\in \R_\alg$ $-$then
\[
\deg(x \pe) \ = \ \deg(\pe).
\]

[Consider $\sqrt[\leftroot{-3}\uproot{3}k]{\abs{x}} \in \ \R_\alg^\times$ \ ($k = \deg(\pe)$).]
\end{x}
\vspace{0.3cm}

Period theory leads to some transcendental conclusions.
\\
\vspace{0.3cm}

%%----------------------------------------------------------------------------------------------06
\begin{x}{\small\bf THEOREM} \ %15
Let $\pe_1$, $\pe_2$ be transcendental periods.  Assume: 
\[
\deg(\pe_1) \ \neq \  \deg(\pe_2).
\]
Then $\pe_1 + \pe_2$ is a transcendental number.

PROOF \ 
Assume false, so
\[
\deg(\pe_1 + \pe_2) \ = \ 0 
\qquad 
\text{or} 
\qquad
\deg(\pe_1 + \pe_2) \ = \ 1.
\]
But
\begin{align*}
\deg(\pe_1 + \pe_2) =  0 \ 
&\implies \ 
\pe_1 + \pe_2 = 0
\\[15pt]
&\implies \ 
\pe_1 = -\pe_2
\\[15pt]
&\implies \ 
\deg(\pe_1) = \deg(\pe_1),
\end{align*}
leaving
\[
\deg(\pe_1 + \pe_2) \ = \ 1,
\]
thus $\pe_1 + \pe_2 \in \ \R_\alg^\times$ (cf. $\# 10$), say $\pe_1 + \pe_2 = x \in \R_\alg^\times$, from which
\begin{align*}
\deg(\pe_1) \ 
&=\ 
\deg(x - \pe_2) 
\\[15pt]
&=\
\deg((-1)\pe_2)\qquad  \text{(cf. $\#13$)}
\\[15pt]
&=\
\deg(\pe_2) \qquad \qquad \text{(cf. $\#14$)}.
\end{align*}
Contradiction.
\end{x}
\vspace{0.3cm}

\begin{x}{\small\bf EXAMPLE} \ %16
It is conjectured that $e + \pi$ is transcendental.  
Recalling that $\deg(\pi) = 2$ (cf. $\# 7$), suppose for sake of argument that 
$e \in \sP_\text{KZ}$ and $\deg(e) \geq 3$, hence $\deg(e) \neq \deg(\pi)$.  
Recalling that $e$ and $\pi$ are transcendental, it follows that $e + \pi$ is transcendental.

[Note: \ 
Nevertheless, the conjecture is that $e$ is not a period.]
\end{x}
\vspace{0.3cm}

\begin{x}{\small\bf THEOREM} \ %17
Let $\pe_1$, $\pe_2$ be transcendental periods.  Assume:
\[
\deg(\pe_1) \ \neq \  \deg(\pe_2).
\]
%%----------------------------------------------------------------------------------------------07
Then $\pe_1 / \pe_2$ is a transcendental number.

[For $\pe_1 / \pe_2 = x \in \ \R_\alg^\times$ implies that
\[
\deg(\pe_1) \ = \ \deg (x \pe_2) \ = \ \deg(\pe_2).]
\]
\end{x}
\vspace{0.3cm}

%%----------------------------------------------------------------------------------------------07
%%%%%%%%%%%%%%%%%%%%%%%%%%%%%%%%%%%%%%
%%%%%%%%%%%%%%%%%%%%%%%%%%%%%%%%%%%%%%
%%%%%%%%%%%%%%%%%%%%%%%%%%%%%%%%%%%%%%

%123
\chapter{
$\boldsymbol{\S}$\textbf{5}.\quad  PRIMITIVE RECURSIVE FUNCTIONS}
\setlength\parindent{2em}
\setcounter{theoremn}{0}
\renewcommand{\thepage}{\S5-\arabic{page}}
%%----------------------------------------------------------------------------------------------01

\begin{x}{\small\bf \un{N.B.}} \ %01
In this $\S$ (and all subsequent ones), $\N$ will stand for $\{0, 1, 2, \ldots \}$ and not $\{1, 2, \ldots \}$.  
Elements of $\N$ will be denoted by $x$, $y$, $z$ or $a$, $b$, $c$ or $n$, $m$, $k$, $\ell$.
\end{x}

\begin{x}{\small\bf \un{N.B.}} \ %02
In what follows, it will be a question of functions $f:\N^n \ra \N$ $(n = 0, 1, 2, \ldots)$.
\end{x}

\begin{x}{\small\bf DEFINITION} \ %03

\hspace{0.5cm} \textbullet \quad 
The \un{zero function} $Z:\N \ra \N$, $Z(x) = 0$.
\\[-.25cm]

\hspace{0.5cm} \textbullet \quad 
The \un{successor function} $S:\N \ra \N$, $S(x) = x + 1$.
\\[-.25cm]

\hspace{0.5cm} \textbullet \quad 
The \un{projection functions} $P_i^n:\N^n \ra \N$, $P_i^n(x_1, \ldots, x_n) = x_i$.
\\[-.25cm]

[Note: \ 
These functions are the so-called \un{initial functions}.]
\end{x}
\vspace{.2cm}

\begin{x}{\small\bf EXAMPLE} \ %04
Constant functions are built up from \mS and \mZ:
\begin{align*}
& S(Z(x))  \hspace{0.7cm} = 1, 
\\[8pt]
&S(S(Z(x))) \ = 2, 
\\[8pt]
&\text{etc.}
\end{align*}
\end{x}
\vspace{0.2cm}

\begin{x}{\small\bf EXAMPLE} \ %05
The \un{addition function}$\text{add}:\N^2 \ra \N$  is defined by
\[
\text{add}(x,y) 
\ = \ 
x + y.
\]
Here
\begin{align*}
& \text{add}(0,y) \hspace{0.8cm} = \ y, 
\\[8pt]
& \text{add}(x + 1,y) \ = \ S(\text{add}(x,y)).
\end{align*}

[Note: \ 
\[
S \quad \implies \quad \text{add}.
\]
For
\[
x + y  
\ = \ 
x + 1 +  1 + \cdots + 1 \hspace{1cm}
\]
\\[-1.1cm]
\begin{tikzcd}%[sep=small]
{\hspace{6.5cm}} \arrow[rrr,dash]
&&&{}
\end{tikzcd}
.]

\hspace{7cm} $y$
\end{x}
\vspace{0.3cm}
%%----------------------------------------------------------------------------------------------02

\begin{x}{\small\bf EXAMPLE} \ %06
The \un{multiplication function} 
$\text{mul}:\N^2 \ra \N$ 
 is defined by
\[
\text{mul}(x,y) 
\ = \ 
x y.
\]
Here
\begin{align*}
& \text{mul}(0,y) = 0, 
\\[8pt]
& \text{mul}(x + 1, y) = \text{add}(\text{mul}(x,y),y).
\end{align*}

[Note: \ 
\[
\text{add \ $\implies$ \ mul.}
\]
For
\[
x y 
\ = \ 
x + x + x + \cdots + x \hspace{1cm}
\]
\\[-1.1cm]
\begin{tikzcd}%[sep=small]
{\hspace{5.6cm}} \arrow[rrrr,dash]
&&&&{}
\end{tikzcd}
.]

\hspace{7cm} $y$
\end{x}
\vspace{0.3cm}

\begin{x}{\small\bf DEFINITION} \ %07
The \un{modified subtraction function} 
$\text{sub}:\N^2 \ra \N$ is the prescription
\[
\text{sub}(x,y) 
\ = \ 
\begin{tikzcd}[sep=small]
{x} \arrow[rr,dash,"\text{\textbullet}"]
&&{y}
\end{tikzcd}
\ = \ 
\begin{cases}
\ x - y \hspace{0.6cm} \text{if} \ x \geq y\\[8pt]
\ \ \ 0 \hspace{1cm} \text{if} \ x < y\
\end{cases}
.
\]
\end{x}
\vspace{0.2cm}

\begin{x}{\small\bf EXAMPLE} \ %08
The \un{absolute value function} $\text{AV}:\N^2 \ra \N$ is defined by
\[
\text{AV} (x,y) 
\ = \ 
\abs{x - y}
\]
and
\[
\abs{x - y} 
\ = \ 
\text{add}
\begin{tikzcd}[sep=small]
{(x} \arrow[rr,dash,"\text{\textbullet}"]
&&{y,}
\end{tikzcd}
\begin{tikzcd}[sep=small]
{y} \arrow[rr,dash,"\text{\textbullet}"]
&&{x)}
\end{tikzcd}
.
\]
\end{x}
\vspace{0.3cm}

\begin{x}{\small\bf EXAMPLE} \ \ %09
$
\text{The} \ 
\begin{cases}
\ \text{minimum}\\
\ \text{maximum}
\end{cases}
\quad
\text{function} 
\quad \ 
\begin{cases}
\ \min:\N^2 \ra \N\\
\ \max:\N^2 \ra \N
\end{cases}
\quad
$

\noindent is defined by
\[
\begin{cases}
\ \min(x,y) \ = \ 
\begin{tikzcd}[sep=small]
{x} \arrow[rr,dash,"\text{\textbullet}"]
&&{}
\end{tikzcd}
\begin{tikzcd}[sep=small]
{(x} \arrow[rr,dash,"\text{\textbullet}"]
&&{y)}
\end{tikzcd}
\ = \ 
\begin{tikzcd}[sep=small]
{y} \arrow[rr,dash,"\text{\textbullet}"]
&&{}
\end{tikzcd}
\begin{tikzcd}[sep=small]
{(y} \arrow[rr,dash,"\text{\textbullet}"]
&&{x)}
\end{tikzcd}
\\[8pt]
\ \max(x,y) \ = \ 
y + 
\begin{tikzcd}[sep=small]
{(x} \arrow[rr,dash,"\text{\textbullet}"]
&&{y)}
\end{tikzcd}
\hspace{1.1cm}
\ = \ 
x + 
\begin{tikzcd}[sep=small]
{(y} \arrow[rr,dash,"\text{\textbullet}"]
&&{x)}
\end{tikzcd}
\end{cases}
.
\]
\end{x}
\vspace{0.3cm}
%%----------------------------------------------------------------------------------------------03

\begin{x}{\small\bf EXAMPLE} \ %10
The \un{identity function} $\text{id}_\N$ can be obtained from \mS and \mZ and sub:
\[
\begin{tikzcd}[sep=small]
{S(x)} \arrow[rr,dash,"\text{\textbullet}"]
&&{S \circ Z(x)}
\end{tikzcd}
.
\]
\end{x}
\vspace{0.3cm}

\begin{x}{\small\bf DEFINITION} \ %11
Given a function $g:\N^m \ra \N$ and given functions $h_i:\N^n \ra \N$ $(i = 1, \ldots, m)$, 
the \un{composition of $g$ and the $h_i$} is the function $f:\N^n \ra \N$, denoted by
\[
g \circ (h_1, \ldots, h_m),
\]
such that 
\[
f(x_1, \ldots, x_n)
\ = \ 
g(h_1(x_1, \ldots, x_n), \ldots, h_m(x_1, \ldots, x_n)).
\]
\end{x}
\vspace{0.3cm}

\begin{x}{\small\bf EXAMPLE} \ %12
Take $m = 2$, $g:\N^2 \ra \N$ and take $n = 2$, $h_1 = P_2^2$, $h_2 = P_1^2$ $-$then
\[
f 
\ = \ 
g \circ \big(P_2^2, P_1^2\big)
\]
is given by
\begin{align*}
f(x_1,x_2) \ 
&=\ 
g(P_2^2(x_1,x_2), P_1^2(x_1,x_2))
\\[11pt]
&=\ 
g(x_2,x_1).
\end{align*}
\end{x}
\vspace{0.3cm}

\begin{x}{\small\bf DEFINITION} \ %13
Suppose that $g:\N^{m-1} \ra \N$, $h:\N^{m+1} \ra \N$ $-$then the function $f:\N^m \ra \N$ is said to be obtained by 
\un{primitive recursion} from $g$ and $h$ if
\[
f(0,x_2, \ldots, x_m) 
\ = \ 
g(x_2, \ldots, x_m)
\]
and 
\[
f(x+1, x_2, \ldots, x_m)
\ = \ 
h(x, f(x, x_2, \ldots, x_m), x_2, \ldots, x_m)
\]
for all $x, x_2, \ldots, x_m \in \ \N$.
\end{x}
\vspace{0.3cm}
%%----------------------------------------------------------------------------------------------04

\begin{x}{\small\bf EXAMPLE} \ %14
Take $m = 1$ $-$then $f:\N \ra \N$ is obtained from the constant $c \in \N$ $(g(0) = c)$ and 
$h:\N^2 \ra \N$ by primitive recursion of $f(0) = c$ and $f(x+1) = h(x,f(x))$.
\end{x}
\vspace{0.3cm}

\begin{x}{\small\bf DEFINITION} \ %15
If $\sF$ is a set of natural number functions and if $\Omega$ is a set of operators on natural number functions, then 
$\text{clos}(\sF,\Omega)$ is the \un{inductive closure} of $\sF$ with respect to $\Omega$, i.e., the smallest set of natural number functions containing 
$\sF$ and closed with respect to the operations of $\Omega$.
\end{x}
\vspace{0.3cm}

\begin{x}{\small\bf DEFINITION} \ %16
The set of \un{primitive recursive functions}, denoted $\fP\fR$, 
is the inductive closure of the initial functions with respect to the operations of composition and primitive recursion.
\end{x}
\vspace{0.3cm}

\begin{x}{\small\bf EXAMPLE} \ %17
The \un{factorial function} is primitive recursive:
\begin{align*}
&0! = 1, 
\\[15pt]
&(n+1)! = \text{mul} \circ (S \circ P_1^2, P_2^2) (n,n!).
\end{align*}
\end{x}
\vspace{0.3cm}

\begin{x}{\small\bf REMARK} \ %18
All the functions encountered in this $\S$ are primitive recursive (but see below).
\end{x}
\vspace{1.cm}

\[
\text{APPENDIX}
\]
\vspace{0.5cm}

Define the \un{Ackermann function} $A:\N^2 \ra \N$ by the following equations:
\begin{align*}
A(0,n) \hspace{1.5cm} 
&= \ 
n+1 
\\[11pt]
A(m+1,0)  \hspace{0.75cm} 
&=  \ 
A(m,1)
\\[11pt]
A(m+1, n+1) \ 
&=  \ 
A(m,A(m+1,n)).
\end{align*}
Then \mA is not primitive recursive.

%%----------------------------------------------------------------------------------------------07
%%%%%%%%%%%%%%%%%%%%%%%%%%%%%%%%%%%%%%
%%%%%%%%%%%%%%%%%%%%%%%%%%%%%%%%%%%%%%
%%%%%%%%%%%%%%%%%%%%%%%%%%%%%%%%%%%%%%

%123
\chapter{
$\boldsymbol{\S}$\textbf{6}.\quad  ELEMENTARY FUNCTIONS}
\setlength\parindent{2em}
\setcounter{theoremn}{0}
\renewcommand{\thepage}{\S6-\arabic{page}}
%%----------------------------------------------------------------------------------------------01

\begin{x}{\small\bf \un{N.B.}} \ %01
All conventions introduced in $\S 5$ remain in force.
\end{x}

\begin{x}{\small\bf DEFINITION}  \ %02
Let $f(t, x_1, \ldots, x_n)$ be a function with $n + 1$ arguments.
\\[-.25cm]

\hspace{0.5cm} \textbullet \quad
\un{bounded summation} is defined by
\[
\sum\limits_{t \leq x} \hsx f(t, x_1, \ldots, x_n) 
\ = \ 
f(0, x_1, \ldots, x_n) + \cdots + f(x, x_1, \ldots, x_n).
\]

\hspace{0.5cm} \textbullet \quad
\un{bounded product} is defined by
\[
\prod\limits_{t \leq x} \hsx f(t, x_1, \ldots, x_n) 
\ = \ 
f(0, x_1, \ldots, x_n) \times \cdots \times f(x, x_1, \ldots, x_n).
\]
\end{x}

\begin{x}{\small\bf DEFINITION} \ %03
The set of \un{elementary functions}, denoted by EL, is the inductive closure of the initial functions, 
the addition, the multiplication, and the modified subtraction with respect to the operations of 
composition, bounded summation, and bounded product.
\end{x}
\vspace{0.2cm}

\begin{x}{\small\bf THEOREM} \ %04
\[
\text{EL} \ \subset \ \fP\fR.
\]

[Note: \ 
The containment is strict.  Consider, e.g., 
\[
f(x) 
\ = \ 
%x^{x^{\reflectbox{$\ddots$}^x}}.]
x^{{\text{\small{$x$}}}^{\cdot^{\cdot^{\cdot^{\text{{\small{$x$}}}}}}}}.]
\]
\end{x}
\vspace{0.2cm}

\begin{x}{\small\bf EXAMPLE} \ %05
\[
\sgn(x) \ = \ 
\begin{cases}
\ 1 \quad \text{if} \quad x \neq 0\\[8pt]
\ 0 \quad \text{if} \quad x = 0
\end{cases}
= 
\begin{tikzcd}[sep=small]
{1} \arrow[rr,dash,"\text{\textbullet}"]
&&{}
\end{tikzcd}
\begin{tikzcd}[sep=small]
{(1} \arrow[rr,dash,"\text{\textbullet}"]
&&{x)}
\end{tikzcd}
\]
is elementary.
\end{x}
\vspace{0.3cm}
%%----------------------------------------------------------------------------------------------02

\begin{x}{\small\bf EXAMPLE} \ %06
\[
> (x,y) \ = \ 
\begin{cases}
\ 1 \quad \text{if} \quad x > y\\[8pt]
\ 0 \quad \text{if} \quad x \leq y
\end{cases}
= 
\begin{tikzcd}[sep=small]
{\sgn(x} \arrow[rr,dash,"\text{\textbullet}"]
&&{y)}
\end{tikzcd}
\]
is elementary.
\\[-0.25cm]

[Ditto for $\geq$, $<$, $\leq$.]
\end{x}
\vspace{0.3cm}

\begin{x}{\small\bf EXAMPLE} \ %07
\[
\text{div}(x,y) 
\ = \ 
\bigg[\frac{x}{y + 1}\bigg]
\]
is elementary.

[In fact, 
\[
\text{div}(x,y) 
\ = \ 
\begin{tikzcd}[sep=small]
{\bigg(\ds\sum\limits_{i=0}^x \geq (x,i(y+1))\bigg)} \arrow[rr,dash,"\text{\textbullet}"]
&&{1}
\end{tikzcd}
.]
\]
\end{x}
\vspace{0.2cm}

\begin{x}{\small\bf EXAMPLE} \ %08
\[
\text{pow}(n,m) \ = \ n^m
\]
is elementary.
\end{x}
\vspace{0.3cm}

\begin{x}{\small\bf NOTATION} \ %09
Given a set $\sF$ of natural number functions, let
\[
\sF^*
\ = \ 
\sF \ \cup \ \big\{P_i^n\big\}
\]
and put
\[
\text{clos}(\sF) 
\ = \ 
\text{clos}(\sF^*, \text{\{composition\}}).
\]

[Note: \ 
Here, in the notation of $\S 5$, $\# 13$, $\Omega$ is the operation of composition.]
\end{x}
\vspace{0.3cm}

\begin{x}{\small\bf THEOREM} \ %10
\[
\text{EL} 
\ = \ 
\text{clos(s, sub, div, pow)}.
\]

[This result is due to Stefano Mazzanti
\footnote[2]{\textit{Mathematical Logic Quarterly} {\bf{48}} (2002), pp. 93-104.}
.]
\end{x}
\vspace{0.3cm}

%%----------------------------------------------------------------------------------------------03
\begin{x}{\small\bf DEFINITION} \ %11
The set of \un{lower elementary functions}, denoted $\ell$EL, is the inductive closure of the initial functions, 
the addition, the multiplication, and the modified subtraction with respect to the operations of composition and bounded summation.
\end{x}
\vspace{0.3cm}

\begin{x}{\small\bf \un{N.B.}} \ %12
Obviously
\[
\ell\text{EL} \ \subset \ \text{EL}.
\]

[Note: \ 
The containment is strict.  Consider, e.g., 
\[
f(n) \ = \ 2^n.]
\]
\end{x}
\vspace{0.3cm}

%%----------------------------------------------------------------------------------------------07
%%%%%%%%%%%%%%%%%%%%%%%%%%%%%%%%%%%%%%
%%%%%%%%%%%%%%%%%%%%%%%%%%%%%%%%%%%%%%
%%%%%%%%%%%%%%%%%%%%%%%%%%%%%%%%%%%%%%

%123
\chapter{
$\boldsymbol{\S}$\textbf{7}.\quad  HIERARCHIES}
\setlength\parindent{2em}
\setcounter{theoremn}{0}
\renewcommand{\thepage}{\S7-\arabic{page}}
%%----------------------------------------------------------------------------------------------01

\ 
\\[-1cm]
\indent

Here is an approach to $\sE^n$ $(n \in \ \N)$, the $n^\text{th}$ set in the Grzegorczyk hierarchy.  
\\[-0.25cm]

[Note: \ Recall that $n$ runs through $0, 1, 2, \ldots$ and not $1, 2, \ldots$ .]
\\

\begin{x}{\small\bf NOTATION} \ %01
\begin{align*}
&f_0(x,y) \ = \ x + 1,\\[3pt]
&f_1(x,y) \ = \ x +y,\\[3pt]
&f_2(x,y) \ = \ x y,\\[3pt]
&f_{n+1}(x,0) \ = \ 1,\\[3pt]
&f_{n+1}(x,y+1) \ = \ f_n(x,f_{n+1}(x,y)) \qquad (n \geq 2).
\end{align*}
\end{x}

\begin{x}{\small\bf DEFINITION}  \ %02
$\sE^n$ is the inductive closure of the initial functions and the $n^\text{th}$ $f_n$ with respect to the operations of composition and bounded primitive recursion.
\end{x}

\begin{x}{\small\bf \un{N.B.}} \ %03
Suppose that $g:\N^{m-1} \ra \N$, $h:\N^{m+1} \ra \N$, $j:\N^{m} \ra \N$,  are in $\sE^n$, 
that $f$ is obtained by primitive recursion from $g$ and $h$, and in addition,
\[
f(x, x_2, \ldots, x_m) 
\ \leq \ 
j(x, x_2, \ldots, x_m).
\]
Then $f$ is in $\sE^n$ as well and $f$ is said to be obtained by  \un{bounded primitive recursion} from $g$, $h$, $j$.
\end{x}
\vspace{0.2cm}

\begin{x}{\small\bf REMARK} \ %04
If $g$ and $h$ belong to $\sE^n$ and $f$ is obtained by primitive recursion from $g$ and $h$, then $f$ belongs to $\sE^{n+1}$.
\end{x}
\vspace{0.2cm}

\begin{x}{\small\bf THEOREM} \ %05
$\forall \ n \in \ \N$, $\sE^n$ is a proper subset of $\sE^{n+1}$.
\end{x}
\vspace{0.3cm}

\begin{x}{\small\bf THEOREM} \ %06
$\forall \ n \in \ \N$, $\sE^n \subset \fP\fR$ and 
\[
\bigcup\limits_n \hsx \sE^n 
\ = \ 
\fP\fR.
\]
\end{x}
\vspace{0.3cm}
%%----------------------------------------------------------------------------------------------02

Therefore the sets
\[
\sE^0, \ \sE^1 - \sE^0, \ \sE^2 - \sE^1, \ldots, \ \sE^n - \sE^{n-1}, \ldots
\]
partition $\fP\fR$.
\\

\begin{x}{\small\bf THEOREM} \ %07
\[
\sE^3 \ = \ \text{EL}.
\]
\end{x}
\vspace{0.2cm}

\begin{x}{\small\bf THEOREM} \ %08
\[
\sE^2 \ \supset \ \ell\text{EL}.
\]

[Note: \ 
It is not known whether all functions in $\sE^2$ are lower elementary.]
\end{x}
\vspace{0.3cm}

\begin{x}{\small\bf \un{N.B.}} \ %09
For the record, $\ell\text{EL}$ is the inductive closure of the inital functions, the addition, the multiplication, and the modified subtraction 
with respect to the operations of composition and bounded summation.
\end{x}
\vspace{0.3cm}

\begin{x}{\small\bf \un{N.B.}} \ %10
For the record, $\sE^2$ is the inductive closure of the inital functions, the addition, the multiplication, and the modified subtraction 
with respect to the operations of composition and bounded primitive recursion.
\end{x}
\vspace{0.3cm}

\begin{x}{\small\bf LEMMA} \ %11
Bounded summation can be derived from bounded primitive recursion.
\end{x}
\vspace{0.3cm}

Consequently 
\[
\sE^2 \ \supset \ \ell\text{EL}.
\]

%%----------------------------------------------------------------------------------------------07
%%%%%%%%%%%%%%%%%%%%%%%%%%%%%%%%%%%%%%
%%%%%%%%%%%%%%%%%%%%%%%%%%%%%%%%%%%%%%
%%%%%%%%%%%%%%%%%%%%%%%%%%%%%%%%%%%%%%

%123
\chapter{
$\boldsymbol{\S}$\textbf{8}.\quad  COMPUTABILITY}
\setlength\parindent{2em}
\setcounter{theoremn}{0}
\renewcommand{\thepage}{\S8-\arabic{page}}
%%----------------------------------------------------------------------------------------------01

\ 
\\[-1cm]
\indent 

Let $\sF$ be a set of natural number functions.
\\

\begin{x}{\small\bf DEFINITION} \ %01
$\sF$ is said to satisfy the \un{standard conditions} if it contains 
the intial functions, the addition, the multiplication, and the modified subtraction and is closed under composition.
\end{x}

\begin{x}{\small\bf \un{N.B.}} \ %02
Both $\ell$EL and $\sE^2$ satisfy the standard conditions.
\end{x}

In what follows, it will be assumed that $\sF$ satisfies the standard conditions.  

[Note: \ 
One consequence of this is that $\sF$ necessarily contains all polynomials with coefficients from $\N$, 
in particular, $\sF$ contains the constant functions.]

\begin{x}{\small\bf DEFINITION} \ %03
An \un{$\sF$-sequence} is a function $A:\N \ra \Q$ that has a representation of the form
\[
A(x) 
\ = \ 
\frac{f(x) - g(x)}{h(x) + 1} \qquad (x = 0, 1, 2, \ldots),
\]
where $f$, $g$, $h: \N \ra \N$ belong to $\sF$.
\end{x}
\vspace{0.2cm}

\begin{x}{\small\bf EXAMPLE} \ %04
Every rational number $q$ gives rise to an $\sF$-sequence $A:\N \ra \Q$, viz. $\forall \ x$, $A(x) = q$.
\\[-0.25cm]

[Suppose that $\ds x = \frac{r}{s}$ and consider the situation when $r \geq 0$, $s \geq 1$.
\\[-0.25cm]

\hspace{0.5cm} \textbullet \quad 
\un{$s = 1$} \quad Set $f(x) = r$, $g(x) = 0$, $h(x) = 0$, hence
\[
\frac{f(x) - g(x) }{h(x) + 1} \ = \ r.
\]

\hspace{0.5cm} \textbullet \quad 
\un{$s > 1$} \quad Set $f(x) = r$, $g(x) = 0$, $h(x) = s - 1$, hence
%%----------------------------------------------------------------------------------------------02
\[
\frac{f(x) - g(x) }{h(x) + 1} \ = \ \frac{r}{s}.]
\]
\end{x}
\vspace{0.2cm}

\begin{x}{\small\bf SUBLEMMA} \ %05
If $f$, $g$ are one-argument elements of $\sF$, then
\[
f + g \in \ \sF, \quad f \cdot g \in \ \sF.
\]

[Consider
\[
\text{add} \circ (f,g) , \quad \text{mul} \circ (f,g),
\]
bearing in mind that $\sF$ is closed under composition.]
\end{x}
\vspace{0.3cm}

\begin{x}{\small\bf LEMMA} \ %06
If $A$, $B:  \N \ra \Q$ are $\sF$-sequences, then so are $A + B$, $A - B$ and $A \cdot B$.
\end{x}
\vspace{0.3cm}

\begin{x}{\small\bf SUBLEMMA} \ %07
If $h$ is a one-argument element of $\sF$, then $\abs{h} \in \ \sF$.
\\[-0.25cm]

[The fact that the modified subtraction function belongs to $\sF$ implies that the absolute value function belongs to $\sF$ 
(cf. $\S5$, $\#8$).]
\end{x}
\vspace{0.2cm}

\begin{x}{\small\bf LEMMA} \ %08
If $A:\N \ra \Q$ is an $\sF$-sequence and if $\forall \ x$, $A(x) \neq 0$, then $\ds\frac{1}{A}$ belongs to $\sF$.

[In fact,
\[
\frac{1}{A(x)} 
\ = \ 
\frac{(h(x) + 1) f(x) - (h(x) + 1)g(x)}{\abs{\abs{f(x) + g(x)}^2 - 1} + 1}.]
\]
\end{x}
\vspace{0.3cm}

\begin{x}{\small\bf DEFINITION} \ %09
A real number $\alpha$ is said to be \un{$\sF$-computable} if there exists an $\sF$-sequence $A:\N \ra \Q$ such that $\forall \ x$,
\[
\abs{A(x) - \alpha} 
\ \leq \ 
\frac{1}{x + 1}.
\]

[Note: \ 
In practice, it can happen that the relation
\[
\abs{A(x) - \alpha} 
\ \leq \ 
\frac{1}{x + 1}
\]
%%----------------------------------------------------------------------------------------------03
is valid only for $x \geq x_0$.  
To remedy this, let
\[
A(0) \ = \ A(x_0), \ 
A(1) = A(x_0), \ldots, A(x_0 - 1) = A(x_0).
\]
Then for $0 \leq x \leq x_0 -1$,
\[
\abs{A(x) - \alpha}
\ = \ 
\abs{A(x_0) - \alpha}
\ \leq \ 
\frac{1}{x_0 + 1} 
\ \leq \ 
\frac{1}{x + 1}.]
\]
\end{x}
\vspace{0.3cm}

\begin{x}{\small\bf NOTATION} \ %10
Denote the set of all $\sF$-computable real numbers by the symbol $\R_\sF$.
\end{x}
\vspace{0.3cm}

\begin{x}{\small\bf \un{N.B.}} \ %11
The constant functions from $\N$ to $\Q$ are $\sF$-computable, hence $\Q \subset \ \R_\sF$ (in particular, $-1 \in \ \R_\sF$).
\end{x}
\vspace{0.3cm}

\begin{x}{\small\bf EXAMPLE} \ %12
Take for $\sF$ the set $\fP\fR$ of primitive recursive functions $-$then the $\sF$-computable real numbers are the primitive 
recursive real numbers.
\\[-0.25cm]

[To arrive at $\R_\sF = \R$, take instead for $\sF$ the set of all natural number functions.]
\end{x}
\vspace{0.3cm}

\begin{x}{\small\bf LEMMA} \ %13
If \mA is an $\sF$-sequence and if $\phi$ is a one-argument element of $\sF$, then the assignment 
$x \ra A(\phi(x))$ is an $\sF$-sequence.
\\[-0.25cm]

[For
\begin{align*}
A(\phi(x)) \
&= \ 
\frac{f(\phi(x)) - g(\phi(x))}{h(\phi(x)) + 1}
\\[11pt]
&=\ 
\frac{(f \circ \phi) (x) - (g \circ \phi) (x)}{(h \circ \phi)(x) + 1}
\end{align*}
and $\sF$ is closed under composition.]
\end{x}
\vspace{0.3cm}

\begin{x}{\small\bf APPLICATION} \ %14
Supposee that \mA is an $\sF$-sequence and $\alpha$ is a real number.  
Assume: \ $x \abs{A(x) - \alpha}$ is bounded $-$then $\alpha \in \ \R_\sF$.
\\[-0.25cm]
%%----------------------------------------------------------------------------------------------04

PROOF \ 
Choose a positive integer $c$ such that $\forall \ x$, 
\[
x \abs{A(x) - \alpha} \ \leq \ c.
\]
Then
\[
\abs{A(c x  + c) - \alpha} \ \leq \ \frac{c}{cx + c} \ = \ \frac{1}{x + 1}.
\]

[Note: \ 
In $\#13$, take $\phi(x) = cx + c$.]
\end{x}
\vspace{0.3cm}

\begin{x}{\small\bf THEOREM} \ %15
$\R_\sF$ is a field.\\

We shall break the proof up into two parts.
\\

\hspace{0.5cm} \un{PART 1} \quad 
Let $\alpha$, $\beta \in \ \R_\sF$ $-$then 
\begin{align*}
\abs{(A(x) + B(x)) - (\alpha + \beta)}\ 
&\leq \ 
\abs{A(x) - \alpha} + \abs{B(x) - \beta}
\\[11pt]
&\leq \ 
\frac{1}{x + 1} + \frac{1}{x + 1} 
\\[11pt]
&=\ 
 \frac{2}{x + 1}
\end{align*}

\qquad $\implies$
\begin{align*}
x \abs{(A(x) + B(x)) - (\alpha + \beta)}\ 
&\leq \ 
\frac{2x}{x + 1} 
\\[11pt]
&= \ 
\frac{2}{1 + \frac{1}{x}} 
\\[11pt]
&\leq \ 
2
\end{align*}

\qquad $\implies$
\[
\alpha + \beta \in \R_\sF
\]
and 
\begin{align*}
\abs{A(x) B(x) - \alpha \beta} \ 
&\leq \ 
\abs{A(x) - \alpha} \abs{B(x)} + \abs{\alpha} \abs{B(x) - \beta}
\\[11pt]
&\leq \ 
\frac{\abs{\beta} + 1 + \abs{\alpha}}{x + 1}
\end{align*}

\qquad $\implies$
\begin{align*}
x \abs{A(x) B(x) - \alpha \beta} \ 
&\leq \ 
\frac{x}{x + 1} (\abs{\beta} + 1 + \abs{\alpha})
\\[11pt]
&\leq \ 
\abs{\beta} + 1 + \abs{\alpha}
\end{align*}

\qquad $\implies$
\[
\alpha \hsx \beta \in \R_\sF.
\]
\\

%%----------------------------------------------------------------------------------------------05
\hspace{0.5cm} \un{PART 2} \quad 
Let $\alpha \neq 0$ and choose $c \in \N$ such that $(c + 1) \abs{\alpha} \geq 2$ $-$then $\forall \ x \geq c$, 
\begin{align*}
\abs{A(x)} \ 
&\geq\ 
\abs{\alpha} - \abs{\alpha - A(x)} 
\\[11pt]
&\geq\ 
\frac{2}{c + 1} - \frac{1}{x + 1} 
\\[11pt]
&\geq\ 
\frac{1}{c + 1}
\end{align*}
\qquad $\implies$
\[
A(x) \ \neq \ 0
\]
and 
\[
\abs{\frac{1}{A(x)} - \frac{1}{\alpha}}
\ = \ 
\frac{\abs{\alpha - A(x)}}{\abs{A(x)}\abs{\alpha}}
\ \leq \ 
\frac{a}{x + 1},
\]
where $a = (c + 1)^2 / 2$.  
Define now a function $C:\N \ra \Q$ by the prescription
\[
C(k) 
\ = \ 
\frac{1}{A(k + c)} \qquad (k + c \geq c \implies A(k + c) \neq 0).
\]
Then \mC is an $\sF$-sequence and 
\[
\abs{C(k) - \frac{1}{\alpha}}
\ \leq \ 
\frac{a}{k + c + 1}
\]
\qquad \qquad $\implies$
\begin{align*}
k \abs{C(k) - \frac{1}{\alpha}}\ 
&\leq \ 
\frac{ka}{k + c + 1}
\\[11pt]
&\leq \ 
\frac{a}{1 + \frac{c}{k} + \frac{1}{k}}
\\[11pt]
&\leq \ 
a
\end{align*} 
\qquad \qquad $\implies$
\[
\frac{1}{\alpha} \in \R_\sF.
\]

In summary: $\R_\sF$ is a field.
\end{x}
\vspace{0.3cm}

\begin{x}{\small\bf DEFINITION} \ %16
Suppose that $f:\N^{k+1} \ra \N$ $-$then the \un{minimizer} $\mu f$ of $f$ is the function
\[
(x_1, \ldots, x_k, x_{k+1}) 
\ra 
\min\{j \in \ \N : f(x_1, \ldots, x_k, j)  =  0 \vee \  j = x_{k+1}\}.
\]
%%----------------------------------------------------------------------------------------------06

[Note: \ 
Spelled out,
\[
\mu f(x_1, \ldots, x_k, x_{k+1})
\]
is the least $j \leq x_{k+1}$ such that 
\[
f(x_1, \ldots, x_k, j) \ = \ 0
\]
if such a $j$ exists, otherwise
\[
\mu f(x_1, \ldots, x_k, x_{k+1})
\ = \ 
x_{k + 1}
\]
if for every $j \leq x_{k+1}$, 
\[
f(x_1, \ldots, x_k, j) \ > \ 0.]
\]
\end{x}
\vspace{0.3cm}

\begin{x}{\small\bf \un{N.B.}} \ %17
To say that $\sF$ is closed under the minimizer operation simply means that 
\[
f \in \ \sF \implies \mu f \in \sF.
\]
\end{x}
\vspace{0.3cm}

\begin{x}{\small\bf THEOREM} \ %18
Suppose that $\sF$ is closed under the minimizer operation.  Let
\[
\alpha_0 \neq 0, \alpha_1, \ldots, \alpha_{k-1}, \alpha_k \in \ \R_\sF.
\]
Then the real roots of the polynomial
\[
P(X) 
\ = \ 
\alpha_0 X^k + \alpha_1 X^{k-1} + \cdots + \alpha_{k-1} X + \alpha_k
\]
belong to $\R_\sF$.
\\[-0.25cm]

PROOF \ 
Let $\zeta$ be a real root of \mP and without loss of generality, assume that 
$P^\prime(\zeta) \neq 0$.  
Choose rational numbers $a$, $b$, $c$, $d$ such that 
$a < \zeta < b$, $0 < c < d$ subject to
\[
c \abs{X - \zeta} 
\ \leq \ 
\abs{P(X)}
\ \leq \ 
d\abs{X - \zeta}
\]
whenever $a \leq X \leq b$.  
Establish the notation:
\[
A_0 \leftrightarrow \alpha_0, \ 
 A_1 \leftrightarrow \alpha_1, 
 \cdots, 
 A_{k-1} \leftrightarrow \alpha_{k-1}, \ 
 A_k \leftrightarrow \alpha_k
\]
%%----------------------------------------------------------------------------------------------07
and introduce the polynomials
\[
P_x(X)  = \ 
A_0(x)X^k + A_1(x)X^{k-1} + \cdots + A_{k - 1}(x)X + A_k(x) \qquad (x = 0, 1, 2, \ldots).
\]
Choose $q \in \Q$:
\[
\abs{P_x(X) - P(X)} 
\ \leq \ 
\frac{q}{x + 1} \qquad (a \leq X \leq b).
\]
Define now a function $A:\N \ra \Q$ via the following procedure.  
Given any $x \in \ \N$, divide $[a,b]$ into $x + 1$ equally long subintervals.  
Let $M_x$ be the set of midpoints of these subintervals, there being at least one $X \in \ M_x$ such that
\[
\abs{P_x(X)} 
\ \leq \
\frac{d(b - a) + 2q}{2(x + 1)}.
\]

Proof: \ 
Choose $X \in \ \M_x$: 
\[
\abs{X - \zeta} 
\ \leq \ 
\frac{b - a}{2(x + 1)}
\]
\qquad \qquad $\implies$
\begin{align*}
\abs{P_x(X)} \ 
&\leq \ 
\abs{P(X)} + \frac{q }{x + 1}
\\[11pt]
&\leq \ 
d \hsx  \frac{(b - a)}{2(x + 1)} + \frac{q}{x + 1}
\\[11pt]
&= \ 
\frac{d(b - a) + 2q}{2(x + 1)}.
\end{align*}
Let $A(x) \in \ M_x$ be the left most element of $M_x$ with the property that 
\[
\abs{P_x(A(x))} 
\ \leq \ 
\frac{d(b - a) + 2q}{2(x + 1)},
\]
so
\begin{align*}
c \abs{A(x) - \zeta} \ 
&\leq \ 
\abs{P(A(x))} 
\\[11pt]
&\leq \ 
\abs{P_x(A(x))}  + \frac{q }{x + 1}
\\[11pt]
&\leq \ 
\frac{d(b - a) + 4q}{2(x + 1)}.
\end{align*}
%%----------------------------------------------------------------------------------------------08
Therefore the product $x \abs{A(x) - \zeta}$ stays bounded as a function of $x$.  
Accordingly, recalling $\#14$, it remains only to show that $A:\N \ra \Q$ is an $\sF$-sequence.  
To this end, note that
\[
A(x) 
\ = \ a + (b - a) \frac{2 \phi(x) + 1}{2x + 2},
\]
where $\phi(x)$ is the smallest $j \in \ \{0, 1, 2, \ldots, x\}$ such that
\[
\abs{P_x\bigg(a + (b - a) \frac{2 j + 1}{2 x + 2}\bigg)}
\ 
\leq \ 
\frac{d(b - a) + 2 q}{2 (x + 1)}
\]
or still, 
\[
\phi(x) = \min\bigg\{j\in \{0, 1, 2, \ldots, x\}: \abs{P_x(a + (b-a) \frac{2j+ 1}{2x + 2}} - \frac{q^\prime}{x + 1} \leq 0\bigg\}, 
\  \bigg(q^\prime = \frac{d(b-a)}{2} + q\bigg).
\]
Since
\[
P_x\bigg(a + (b - a) \frac{2j + 1}{2x + 2}\bigg) 
\ = \ 
\sum\limits_{i=0}^k \hsx A_i(x) \bigg(a + (b - a) \frac{2j + 1}{2x + 2}\bigg)^{k-i},
\]
the function 
\[
x \ra P_x\bigg(a + (b - a) \frac{2j + 1}{2x + 2}\bigg)
\]
defines an $\sF$-sequence for each $j$, hence can be represented in the form
\[
x \ra \frac{f(x,j) - g(x,j)}{h(x,j) + 1}.
\]
Put
\[
\Phi(x,j) 
\ = \ 
\begin{tikzcd}[sep=small]
{f(x,j)} \arrow[rr,dash,"\text{\textbullet}"]
&&{g(x,j),}
\end{tikzcd}
\]
thus
\[
\phi(x) 
\ = \ \min\{j \in \ \N : \Phi(x,j) = 0 \vee \  j = x\},
\]
so $\phi \in \sF$ and this implies that $A:\N \ra \Q$ is an $\sF$-sequence.
\end{x}
\vspace{0.3cm}
%%----------------------------------------------------------------------------------------------09

\begin{x}{\small\bf APPLICATION} \ %19
\[
\R_\alg 
\ \subset \ \R_\sF \ \ (\subset \ \R).
\]

[In other words, every real algebraic number is an $\sF$-computable real number.]
\end{x}
\vspace{0.3cm}

\begin{x}{\small\bf REMARK} \ %20
It is a fact that the minimizer operation can be derived from bounded summation, hence $\ell\text{EL}$ and $\sE^2$ 
are closed under the minimizer operation (cf. $\#2$).
\end{x}
\vspace{0.3cm}

\begin{x}{\small\bf THEOREM} \ %21
Suppose that $\sF$ is closed under the minimizer operation $-$then $\R_\sF$ is a real closed field.
\\[-0.25cm]

PROOF \ 
$\R_\sF$ is an ordered field.  And:
\\[-0.25cm]

\hspace{0.5cm} \textbullet \quad
Every polynomial of odd degree with coefficients in $\R_\sF$ has at least one root in $\R_\sF$.
\\[-0.25cm]

[Since all data is real, on general grounds such a polynomial has at least one real root $\zeta$. And, in view of $\#18$, 
$\zeta \in \ \R_\sF$.]
\\[-0.25cm]

\hspace{0.5cm} \textbullet \quad
If $\alpha > 0$ is an element of $\R_\sF$, then $\alpha$ has a square root $\sqrt{\alpha}$ in $\R_\sF$.
\\[-0.25cm]

[As a positive real number, $\sqrt{\alpha}$ is a root of the polynomial $X^2 - \alpha$.  
But the coefficients of this polynomial are in $\R_\sF$, hence by $\#18$, $\sqrt{\alpha}$ belongs to $\R_\sF$.]
\end{x}
\vspace{0.3cm}

%%%%%%%%%%%%%%%%%%%%%%%%%%%%%%%%%%%%%%
%%%%%%%%%%%%%%%%%%%%%%%%%%%%%%%%%%%%%%
%%%%%%%%%%%%%%%%%%%%%%%%%%%%%%%%%%%%%%

%123
\chapter{
$\boldsymbol{\S}$\textbf{9}.\quad  THE SKORDEV CRITERION$^\dagger$
}
%\footnote[{\textcolor{white}0}]
%\footnote[0]
\footnotetext{\textit{${ }^\dagger$Journal of Universal Computer Science} {\bf{14}} (2008), pp. 861-875.}
\setlength\parindent{2em}
\setcounter{theoremn}{0}
\renewcommand{\thepage}{\S9-\arabic{page}}

\ 
\\[-1.25cm]
\indent 
%%----------------------------------------------------------------------------------------------01

Given a set $\sF$ of natural number functions, assume as in $\S 8$ that $\sF$ satisfies the standard conditions.
\\

\begin{x}{\small\bf DEFINITION} \ %01
An \un{$\sF$-2-sequence} is a function $A:\N^2 \ra \Q$ that has a representation of the form 
\[
A(x,n) 
\ = \ 
\frac{f(x,n) - g(x,n)}{h(x,n) + 1} \qquad (x, n = 0, 1, 2, \ldots),
\]
where $f$, $g$, $h:\N^2 \ra \N$ belong to $\sF$.
\end{x}

\begin{x}{\small\bf DEFINITION} \ %02
A real valued function $\alpha:\N \ra \R$  is said to be  \un{$\sF$-computable} if there exists an $\sF$-2-sequence 
$A:\N^2 \ra \Q$ such that $\forall \ x$, $\forall \ n$, 
\[
\abs{A(x,n) - \alpha(n)}
\ \leq \ 
\frac{1}{x + 1}.
\]
\end{x}

\begin{x}{\small\bf \un{N.B.}} \ %03
It is clear that $\forall \ n \in \ \N$, the real number $\alpha(n)$ is $\sF$-computable (cf. $\S 8$, $\#9$).  
On the other hand, if $\forall \ n \in \ \N$, $\alpha(n)$ is $\sF$-computable, then there exists an $\sF$-sequence $A_n(x)$ such that 
$\forall \ x$,
\[
\abs{A_n(x) - \alpha(n)}
\ \leq \ 
\frac{1}{x + 1},
\]
so setting 
\[
A(x,n) \ = \ A_n(x)
\]
leads to the conclusion that $\alpha$ is $\sF$-computable.
\end{x}
\vspace{0.2cm}

\begin{x}{\small\bf METHODOLOGY} \ %04
Given an $\sF$-sequence $A:\N \ra \Q$ , view it as a function $\alpha:\N \ra \R$ $-$then $\alpha$ is $\sF$-computable, 
i.e., there exists an $\sF$-2-sequence $A:\N^2 \ra \Q$ 
%%----------------------------------------------------------------------------------------------02
such that 
\[
\abs{A(x,n) - \alpha(n)}
\ \leq \ 
\frac{1}{x + 1}.
\]
Thus let
\[
A(x,n) \ = \ A(n)
\]
to get
\begin{align*}
\abs{A(x,n) - \alpha(n)} \ 
&=\ 
\abs{A(n) - A(n)}
\\[15pt]
&=\ 
0
\\[15pt]
&\leq\ 
\frac{1}{x + 1}.
\end{align*}
\end{x}
\vspace{0.2cm}

\begin{x}{\small\bf LEMMA} \ %05
If $\alpha:\N \ra \R$  is $\sF$-computable, then there exist functions $f:\N^2 \ra \N$ and $g:\N^2 \ra \N$ in $\sF$ 
such that $\forall \ x$, $\forall \ n$,
\[
\abs{\frac{f(x,n) - g(x,n)}{x + 1} - \alpha(n)} 
\ \leq \ 
\frac{1}{x + 1}.
\]

PROOF \ 
Changing notation, start with a relation
\[
\abs{\frac{u(x,n) - v(x,n)}{w(x,n) + 1} - \alpha(n)} 
\ \leq \ 
\frac{1}{x + 1}
\]
per $\#2$.  Introduce
\[
\begin{cases}
\ f_0(x,n) \ = \ u(2x + 1, n) \\[8pt]
\ g_0(x,n) \ = \ v(2x + 1, n) \\[8pt]
\ h_0(x,n) \ = \ w(2x + 1, n) 
\end{cases}
.
\]
Then
\begin{align*}
\abs{\frac{f_0(x,n) - g_0(x,n)}{h_0(x,n) + 1} - \alpha(n)} \ 
&=\ 
\abs{\frac{u(2x + 1, n) - v(2x + 1, n)}{w(2x + 1, n) + 1} - \alpha(n)}
\\[15pt]
%%----------------------------------------------------------------------------------------------03
&\leq\
\frac{1}{2x + 1 + 1}
\\[15pt]
&=\
\frac{1}{2(x + 1)}.
\end{align*}

Define $C:\N^2 \ra \N$ by the rule
\[
C(i,j) 
\ = \ 
\bigg[\frac{i}{j+1} + \frac{1}{2}\bigg],
\]
thus $C \in \ \sF$ and 
\[
\abs{C(i,j) - \frac{i}{j+1}} 
\ \leq \ 
\frac{1}{2}.
\]
Put now
\[
f(x,n) 
\ = \ 
C((x + 1) 
\begin{tikzcd}[sep=small]
{(f_0(x,n)} \arrow[rr,dash,"\text{\textbullet}"]
&&{g_0(x,n)),h_0(x,n))}
\end{tikzcd}
\]
and 
\[
g(x,n) 
\ = \ 
C((x + 1) 
\begin{tikzcd}[sep=small]
{(g_0(x,n)} \arrow[rr,dash,"\text{\textbullet}"]
&&{f_0(x,n)),h_0(x,n))}
\end{tikzcd}
.
\]
Then $f$, $g \in \sF$ and (details below)
\[
\abs{f(x,n) - g(x,n) - (x + 1)\frac{f_0(x,n) - g_0(x,n)}{h_0(x,n) + 1}}
\ \leq \ 
\frac{1}{2}.
\]
Multiply this by $\ds\frac{1}{x+1}$ to get 
\[
\abs{\frac{f(x,n) - g(x,n)}{x + 1} - \frac{f_0(x,n) - g_0(x,n)}{h_0(x,n) + 1}}
\ \leq \ 
\frac{1}{2(x + 1)}.
\]
Therefore
\begin{align*}
\abs{\frac{f(x,n) - g(x,n)}{x + 1} - \alpha(n)}
&=\ 
\bigg|{\frac{f(x,n) - g(x,n)}{x + 1} - \frac{f_0(x,n) - g_0(x,n)}{h_0(x,n) + 1}}
\\[15pt]
& \hspace{4cm} +\frac{f_0(x,n) - g_0(x,n)}{h_0(x,n) + 1} - \alpha(n)\bigg|
\\[15pt]
&\leq\ 
\abs{\frac{f(x,n) - g(x,n)}{x + 1} - \frac{f_0(x,n) - g_0(x,n)}{h_0(x,n) + 1}}
\\[15pt]
& \hspace{4cm} +\abs{\frac{f_0(x,n) - g_0(x,n)}{h_0(x,n) + 1} - \alpha(n)}
\\[15pt]
&\leq\ 
\frac{1}{2(x + 1)} + \frac{1}{2(x + 1)} 
\\[15pt]
&=\ 
\frac{1}{x + 1}.
\end{align*}
%%----------------------------------------------------------------------------------------------04
\end{x}
\vspace{0.3cm}

\begin{x}{\small\bf DETAILS} \ %06
The claim is that
\[
\abs{f(x,n) - g(x,n) - (x + 1) \hsx  \frac{f_0(x,n) - g_0(x,n)}{h_0(x,n) + 1}} 
\ \leq \ 
\frac{1}{2}.
\]

\hspace{0.5cm} \textbullet \quad 
Suppose that
\[
f_0(x,n) \ < \  g_0(x,n).
\]
By definition, 
\[
\begin{tikzcd}[sep=small]
{f_0(x,n)} \arrow[rr,dash,"\text{\textbullet}"]
&&{g_0(x,n)}
\end{tikzcd}
=\ 
\begin{cases}
\ f_0(x,n) - g_0(x,n) \quad \text{if} \quad f_0(x,n) \geq g_0(x,n) \\[8pt]
\  \ \ \  0 \hspace{2.85cm} \text{if} \quad f_0(x,n) < g_0(x,n)
\end{cases}
.
\]
Accordingly
\[
\begin{tikzcd}[sep=small]
{f_0(x,n)} \arrow[rr,dash,"\text{\textbullet}"]
&&{g_0(x,n)}
\ = \ 0.
\end{tikzcd}
\]
Therefore
\begin{align*}
f(x,n) \ 
&=\ 
C(0,h_0(x,n))
\\[12pt]
&=\ 
\bigg[\frac{0}{h_0(x,n) + 1} + \frac{1}{2}\bigg]
\\[12pt]
&=\ 
\bigg[\frac{1}{2}\bigg]
\\[12pt]
&=\ 
0.
\end{align*}
So consider
\[
\abs{-g(x,n) - (x +1) \frac{f_0(x,n) - g_0(x,n)}{h_0(x,n) + 1}} 
\]
%%----------------------------------------------------------------------------------------------05
or still, 
\[
\abs{g(x,n) - (x +1) \frac{g_0(x,n) - f_0(x,n)}{h_0(x,n) + 1}} 
\]
By definition,
\[
g(x,n) 
\ = \ 
C((x + 1) 
\begin{tikzcd}[sep=small]
{(g_0(x,n)} \arrow[rr,dash,"\text{\textbullet}"]
&&{f_0(x,n)),h_0(x,n))}
\end{tikzcd}
.
\]
But here
\[
\begin{tikzcd}[sep=small]
{g_0(x,n)} \arrow[rr,dash,"\text{\textbullet}"]
&&{f_0(x,n)}
\end{tikzcd}
\ = \ 
g_0(x,n) - f_0(x,n)
\]
since
\[
f_0(x,n) \ < \  g_0(x,n).
\]
Recalling that
\[
\abs{C(i,j) - \frac{i}{j + 1}} 
\ \leq \ 
\frac{1}{2},
\]
specialize and take
\[
\begin{cases}
\ i = (x + 1) (g_0(x,n) - f_0(x,n)) \\[8pt]
\ j = h_0(x,n)
\end{cases}
.
\]
Then
\begin{align*}
&\abs{g(x,n) - (x +1) \frac{g_0(x,n) - f_0(x,n)}{h_0(x,n) + 1}} 
\\[15pt]
& \hspace{2cm} =\ 
\abs{C((x + 1) (g_0(x,n) - f_0(x,n)), h_0(x,n)) - (x +1) \frac{g_0(x,n) - f_0(x,n)}{h_0(x,n) + 1}}
\\[15pt]
& \hspace{2cm} \leq\ 
\frac{1}{2}.
\end{align*}

\hspace{0.5cm} \textbullet \quad 
Suppose that
\[
f_0(x,n) \geq g_0(x,n).
\]
%%----------------------------------------------------------------------------------------------06
By definition, 
\[
\begin{tikzcd}[sep=small]
{g_0(x,n)} \arrow[rr,dash,"\text{\textbullet}"]
&&{f_0(x,n)}
\end{tikzcd}
\ = \ 
\begin{cases}
\ g_0(x,n) - f_0(x,n) \quad \text{if} \quad g_0(x,n) \geq f_0(x,n)\\[8pt]
\ 0 \hspace{3.5cm} \text{if} \quad g_0(x,n) < f_0(x,n)
\end{cases}
.
\]
Accordingly
\[
\begin{cases}
\ g_0(x,n)  < f_0(x,n)  \ \implies \
\begin{tikzcd}[sep=small]
{g_0(x,n)} \arrow[rr,dash,"\text{\textbullet}"]
&&{f_0(x,n)}
\end{tikzcd}
\ = \ 0
\\[8pt]
\ g_0(x,n)  = f_0(x,n)  \ \implies \
\begin{tikzcd}[sep=small]
{g_0(x,n)} \arrow[rr,dash,"\text{\textbullet}"]
&&{f_0(x,n)}
\end{tikzcd}
\\[8pt]
\hspace{5.05cm} = \ g_0(x,n) - f_0(x,n) \ = \ 0
\end{cases}
.
\]
Therefore
\begin{align*}
g(x,n) \ 
&=\ 
C(0,h_0(x,n+1))
\\[12pt]
&=\ 
\bigg[\frac{0}{h_0(x,n) + 1} + \frac{1}{2}\bigg]
\\[12pt]
&=\ 
\bigg[\frac{1}{2}\bigg]
\\[12pt]
&=\ 
0.
\end{align*}
So consider
\[
\abs{f(x,n) - (x +1) \frac{f_0(x,n) - g_0(x,n)}{h_0(x,n) + 1}}.
\]
By definition,
\[
f(x,n) 
\ = \ 
C((x + 1) 
\begin{tikzcd}[sep=small]
{(f_0(x,n)} \arrow[rr,dash,"\text{\textbullet}"]
&&{g_0(x,n)),h_0(x,n))}
\end{tikzcd}
.
\]
But here
\[
\begin{tikzcd}[sep=small]
{f_0(x,n)} \arrow[rr,dash,"\text{\textbullet}"]
&&{g_0(x,n)}
\end{tikzcd}
\ = \ 
f_0(x,n) - g_0(x,n)
\]
since
\[
f_0(x,n) \ \geq \  g_0(x,n).
\]
%%----------------------------------------------------------------------------------------------07
Recalling that
\[
\abs{C(i,j) - \frac{i}{j + 1}} 
\ \leq \ 
\frac{1}{2},
\]
specialize and take
\[
\begin{cases}
\ i = (x + 1) (f_0(x,n) - g_0(x,n)) \\[8pt]
\ j = h_0(x,n)
\end{cases}
.
\]
Then
\begin{align*}
\hspace{0.8cm}&\abs{f(x,n) - (x +1) \frac{f_0(x,n) - g_0(x,n)}{h_0(x,n) + 1}} 
\\[15pt]
& \hspace{2cm} =\ 
\abs{C((x + 1) (f_0(x,n) - g_0(x,n)), h_0(x,n)) - (x +1) \frac{f_0(x,n) - g_0(x,n)}{h_0(x,n) + 1}}
\\[15pt]
& \hspace{2cm} \leq\ 
\frac{1}{2}.
\end{align*}
\end{x}
\vspace{0.3cm}

\begin{x}{\small\bf \un{N.B.}} \ %07
The upshot is that in the definition of $\sF$-computability, one can take $h(x,n) = x$.
\end{x}
\vspace{0.2cm}

Thus far the only conditions imposed on $\sF$ are the standard ones but to proceed to the main result it will be assumed henceforth 
that $\sF$ is closed under bounded summation (which, of course, is the case if $\sF = \ell$EL or $\sE^2$).
\\

\begin{x}{\small\bf THEOREM} \ %08
Let $\alpha:\N \ra \R$ be $\sF$-computable, assume that the series 
$\ds\sum\limits_{n=0}^\infty \hsx \alpha(n)$ is convergent, and let $\Xi$ be its sum.  
Suppose that there exists a function $\xi:\N \ra \N$ in $\sF$ such that $\forall \ x \in \ \N$,
\[
\abs{\sum\limits_{n=\xi (x) + 1}^\infty \hsx \alpha(n)} 
\ \leq \ 
\frac{1}{x + 1}.
\]
Then $\Xi$ is $\sF$-computable.
\end{x}
\vspace{0.3cm}

%%----------------------------------------------------------------------------------------------08

\begin{x}{\small\bf LEMMA} \ %09
Let $\alpha:\N \ra \R$ be $\sF$-computable.  Define $\alpha^\Sigma:\N \ra \R$ by setting 
\[
\alpha^\Sigma(m)
\ = \ 
\sum\limits_{n=0}^m \hsx \alpha(n).
\]
The $\alpha^\Sigma$ is $\sF$-computable.
\\[-0.25cm]

PROOF \ 
Per $\# 5$, write
\[
\abs{\frac{f(x,n) - g(x,n)}{x + 1} - \alpha(n)} 
\ \leq \ 
\frac{1}{x + 1}
\]
and define functions
\[
\begin{cases}
\ f^\Sigma:\N^2 \ra \N\\[8pt]
\ g^\Sigma:\N^2 \ra \N
\end{cases}
\]
in $\sF$ by stipulating that 
\[
\begin{cases}
\ f^\Sigma(x,m) \ = \ \sum\limits_{n=0}^m \hsx f(xm + x + m,n)\\[11pt]
\ g^\Sigma(x,m) \ = \ \sum\limits_{n=0}^m \hsx g(xm + x + m,n)
\end{cases}
.
\]
Then
\begin{align*}
\abs{\frac{f(xm + x + m,n) - g(xm + x+m,n)}{xm + x + m + 1} - \alpha(n)} \ 
&\leq \ 
\frac{1}{xm + x + m + 1}
\\[15pt]
\implies \hspace{6.7cm}
&
\\[15pt]
\abs{\frac{f^\Sigma(x,m) - g^\Sigma(x,m)}{xm + x + m + 1} - \alpha^\Sigma(m)} \ 
&\leq \ 
\frac{1}{x + 1}.
\end{align*}
I.e.:
\[
\abs{\frac{f^\Sigma(x,m) - g^\Sigma(x,m)}{h(x,m) + 1} - \alpha^\Sigma(m)} \ 
\ \leq \ 
\frac{1}{x + 1}
\]
if
\[
h(x,m) \ = \ xm + x + m,
\]
thus $\alpha^\Sigma$ is $\sF$-computable.
\\
%%----------------------------------------------------------------------------------------------09

Turning now to the proof of the theorem, per $\# 9$, determine functions
\[
F, \ G, \ H:\N^2 \ra \N
\]
in $\sF$ such that $\forall \ x$, $\forall \ n$, 
\[
\abs{\frac{F(x,n) - G(x,n)}{H(x,n) + 1} - \alpha^\Sigma(n)} \ 
\ \leq \ 
\frac{1}{x + 1}.
\]
Introduce
\[
\begin{cases}
\ u(x) \ = \ F(2x + 1, \xi(2x + 1))\\[8pt]
\ v(x) \ = \ G(2x + 1, \xi(2x + 1))\\[8pt]
\ w(x) \ = \ H(2x + 1, \xi(2x + 1))
\end{cases}
.
\]

Then $u$, $v$, $w \in \ \sF$ and
\begin{align*}
\abs{\frac{u(x) - v(x)}{w(x) + 1} - \Xi} \ 
&=\ 
\abs{\frac{u(x) - v(x)}{w(x) + 1} - \sum\limits_{n=0}^\infty \alpha(n)} 
\\[15pt]
&=\ 
\abs{\frac{u(x) - v(x)}{w(x) + 1} - \sum\limits_{n=0}^{\xi(2x+1)} \alpha(n) 
- \sum\limits_{n=\xi(2x+1) + 1}^\infty \alpha(n)} 
\\[15pt]
&=\ 
\bigg|\frac{F(2x + 1, \xi(2x + 1)) - G(2x + 1, \xi(2x + 1))}{H(2x + 1, \xi(2x + 1))} - \alpha^\Sigma(\xi(2x+1))
\\[15pt]
 &\hspace{3cm}  \hspace{3cm} - \sum\limits_{n = \xi(2x + 1) + 1}^\infty \hsx \alpha(n)\bigg|
\\[15pt]
&\leq\ 
\abs{\frac{F(2x + 1), \xi(2x + 1)) - G(2x + 1, \xi(2x + 1))}{H(2x + 1, \xi(2x + 1))} - \alpha^\Sigma(\xi(2x+1))} 
\\[15pt]
&\hspace{6cm} + \abs{\sum\limits_{n = \xi(2x + 1) + 1}^\infty \hsx \alpha(n)}
\\[15pt]
&\leq\ 
\frac{1}{2x + 1 + 1} + \abs{\sum\limits_{n = \xi(2x + 1) + 1}^\infty \hsx \alpha(n)}
\\[15pt]
&\leq\ 
\frac{1}{2x + 1 + 1} + \frac{1}{2x + 1 + 1}
\\[15pt]
&=\ 
\frac{1}{2x + 2} + \frac{1}{2x + 2}
\\[15pt]
&=\ 
\frac{2}{2x + 2}
\\[15pt]
&=\ 
\frac{1}{x + 1}.
\end{align*}
%%----------------------------------------------------------------------------------------------10
\end{x}
\vspace{0.3cm}

\begin{x}{\small\bf \un{N.B.}} \ %10
Assuming that the series $\ds\sum\limits_{n=0}^\infty \hsx \alpha(n)$ is convergent (in practice, this is invariably a non-issue), there are two points.
\\[-0.25cm]

1. \quad Establish that $\alpha$ is $\sF$-computable (or that $\forall \ n$, $\alpha(n)$ is $\sF$-computable).
\\[-0.25cm]

2. \quad Find $\xi$ and deal with the speed of convergence.
\end{x}
\vspace{0.3cm}

\begin{x}{\small\bf EXAMPLE} \ %11
$e$ is $\sE^2$-computable, i.e., $e \in \ \R_{\sE^2}$.  Thus write
\[
e 
\ = \ 
\sum\limits_{n=0}^\infty \hsx \frac{1}{n!}
\]
and let \ $\ds\alpha(n) = \frac{1}{n!}$, \ hence \ $\alpha:\N \ra \Q \subset \ \R$, \ the claim being that $\alpha$ is $\sE^2$-computable.  \ 
To see this, let
\[
f(x,n) 
\ = \ 
\bigg[ \frac{x}{n!}\bigg].
\]
Then $f \in \ \sE^2$:
\begin{align*}
&f(x,0) \ =\ x, 
\\[5pt]
&f(x,n + 1) \ = \ \bigg[\frac{f(x,n)}{n + 1}\bigg], 
\\[5pt]
&f(x,n) \ \leq x
\end{align*}
and one can quote bounded recursion.  
In addition
\begin{align*}
\abs{\frac{f(x + 1,n)}{x + 1} - \alpha(n)} \ 
&=\ 
\abs{\frac{f(x + 1,n)}{x + 1} - \frac{1}{n!}}
\\[15pt]
&=\ 
\frac{1}{x + 1} \abs{f(x + 1,n) - \frac{x + 1}{n!}}
\\[15pt]
&=\ 
\frac{1}{x + 1} \abs{\bigg[\frac{x + 1}{n!}\bigg] - {\frac{x + 1}{n!}}}
\\[15pt]
&\leq \ 
\frac{1}{x + 1}.
\end{align*}

%%----------------------------------------------------------------------------------------------11
Therefore $\alpha$ is $\sF$-computable.   
It remains to define $\xi:\N \ra \N$ and consider
\[
\sum\limits_{n = \xi(x) + 1}^\infty \hsx \alpha(n) 
\ = \ 
\sum\limits_{n = \xi(x) + 1}^\infty \hsx \frac{1}{n!}.
\]
So put
\begin{align*}
&\xi(0) \ = \ 1, \\
&\xi(1) \ = \ 2, \\
&\xi(x) \ = \ x \qquad (x \geq 2).
\end{align*}
Then
\[
(x = 0) \qquad 
\sum\limits_{n = \xi(0) + 1}^\infty \hsx \frac{1}{n!} 
\ = \ 
\sum\limits_{n = 2}^\infty \hsx \frac{1}{n!} 
\ = \ 
e - 2 
\ \leq \ 
1 
\ = \ 
\frac{1}{0 + 1} \hspace{2.75cm}
\]
\[
(x = 1) \qquad 
\sum\limits_{n = \xi(1) + 1}^\infty \hsx \frac{1}{n!} 
\ = \ 
\sum\limits_{n = 3}^\infty \hsx \frac{1}{n!} 
\ = \ 
e - 2.5
\ \sim \ 
2.7 - 2.5
\ = \ 
.2
\ \leq \ 
\frac{1}{1 + 1}
\]
\begin{align*}
(x > 1) \qquad  \sum\limits_{n = \xi(x) + 1}^\infty \hsx \frac{1}{n!} \ 
&=\ 
\sum\limits_{n = x + 1}^\infty \hsx \frac{1}{n!} \hspace{6.4cm}
\\[15pt]
&=\ 
\frac{1}{(x + 1)!} + \frac{1}{(x + 2)!} + \cdots
\\[15pt]
&=\ 
\frac{1}{x!}\bigg( \frac{1}{x + 1} + \frac{1}{(x + 1)(x + 2)} + \cdots
\\[15pt]
&<\ 
\frac{1}{x!}\bigg(\frac{1}{x + 1} + \frac{1}{(x + 1)^2} + \cdots
\\[15pt]
&=\ 
\frac{1}{x!}\bigg(\frac{\frac{1}{x + 1}}{1 - \frac{1}{x + 1}}\bigg)
\\[15pt]
&=\ 
\frac{1}{x!}\frac{1}{x}
\\[15pt]
&\leq\ 
\frac{1}{x + 1}.
\end{align*}
%%----------------------------------------------------------------------------------------------12
Conclusion:
\[
\Xi \ = \ e
\]
is $\sE^2$-computable.  

[Note: \ 
It turns out that $e$ is actually $\ell$EL-computable (cf. $\S 10$, $\# 6$).]
\end{x}
\vspace{0.3cm}

\begin{x}{\small\bf REMARK} \ %12
Recall that $\R_\alg \subset \ R_{\sE^2}$ (cf. $\S 8$, $\#19$) and since $e$ is transcendental, it follows that the containment is proper.
\end{x}
\vspace{0.3cm}

%%%%%%%%%%%%%%%%%%%%%%%%%%%%%%%%%%%%%%
%%%%%%%%%%%%%%%%%%%%%%%%%%%%%%%%%%%%%%
%%%%%%%%%%%%%%%%%%%%%%%%%%%%%%%%%%%%%%

%123
\chapter{
$\boldsymbol{\S}$\textbf{10}.\quad  TECHNICALITIES}
\setlength\parindent{2em}
\setcounter{theoremn}{0}
\renewcommand{\thepage}{\S10-\arabic{page}}
%%----------------------------------------------------------------------------------------------01

\begin{x}{\small\bf DEFINITION} \ %01
A relation $\R \subset \N^n$ is said to be \un{lower elementary} if its characteristic function belongs to $\ell$EL.
\end{x}

\begin{x}{\small\bf LEMMA}  \ %02
Suppose that $f:\N \ra \N$ is in $\ell$EL.  
Define a function $\phi:\N \ra \N$ by the formula
\[
\phi(n) 
\ = \ 
\prod\limits_{k=0}^n \hsx f(k).
\]
Then the graph of $\phi$ is lower elementary.

[Note: \ 
Recall that  $\ell$EL is closed under bounded summation.]
\end{x}

\begin{x}{\small\bf EXAMPLE} \ %03
Fix a positive natural number $N$ and define $f:\N \ra \N$ by stipulating that $\forall \ k$, $f(k) = N$ $-$then $f \in  \ \ell$EL.  
Moreover $\phi(n) = N^{n+1}$ and the graph of $\phi$ is lower elementary.
\end{x}
\vspace{0.2cm}

\begin{x}{\small\bf LEMMA} \ %04
Suppose that $\phi:\N \ra \N$ has the property that $\forall \ n$, $\phi(n) \neq 0$.  
Assume further that the graph of $\phi$ is lower elementary $-$then the function 
\[
n \ra \frac{1}{\phi(n)}
\]
is $\ell$EL-computable.
\end{x}
\vspace{0.2cm}

\begin{x}{\small\bf EXAMPLE} \ %05
For every positive natural number \mN, the function
\[
n \ra \frac{1}{N^{n+1}}
\]
is $\ell$EL-computable.
\end{x}
\vspace{0.3cm}
%%----------------------------------------------------------------------------------------------02

\begin{x}{\small\bf REMARK} \ %06
It was shown in $\S9$, $\#11$ that
\[
e \ = \ \sum\limits_{n=0}^\infty \hsx \frac{1}{n!}
\]
is $\sE^2$-computable.  However more is true: $e$ is $\ell$EL-computable.  
To see this, consider $f:\N \ra \N$, where
\[
f(0) = 1, \quad f(k) = k \qquad (k > 0).
\]
Then
\begin{align*}
\phi(n) \ 
&=\ 
\prod\limits_{k=0}^n \hsx f(k)
\\[11pt]
&=\ 
 f(0) f(1) f(2) \cdots f(n)
\\[11pt]
&=\ 
n!.
\end{align*}
Therefore the function
\[
n \ra \frac{1}{n!}
\]
is $\ell$EL-computable (and the argument proceeds \ldots).
\\[-0.25cm]

[Note: \ 
It is clear that $f$ is in is $\ell$EL.]
\end{x}
\vspace{0.1cm}

\begin{x}{\small\bf LEMMA} \ %07
If $\alpha:\N \ra \R$, $\beta:\N \ra \R$ are $\ell$EL-computable bounded functions, then the product 
$\alpha \hsx \beta :\N \ra \R$ is also $\ell$EL-computable.
\end{x}
\vspace{0.2cm}

\begin{x}{\small\bf EXAMPLE} \ %08
The function 
\[
n \ra \frac{(-1)^n}{N^{n+1}}
\]
is $\ell$EL-computable.
\\[-0.25cm]

[Note: \ 
\begin{align*}
(-1)^n \ 
&=\ 
(n + 1)  \hsx \modx 2 - n \hsx \modx 2
\\[11pt]
&=\ 
\text{``$f(n) - g(n)$''},
\end{align*}
hence $(-1)^n$ is an $\ell$EL-sequence, hence is $\ell$EL-computable (cf. $\S9$, $\#4$).]
\end{x}
\vspace{0.3cm}

%%%%%%%%%%%%%%%%%%%%%%%%%%%%%%%%%%%%%%
%%%%%%%%%%%%%%%%%%%%%%%%%%%%%%%%%%%%%%
%%%%%%%%%%%%%%%%%%%%%%%%%%%%%%%%%%%%%%

%123
\chapter{
$\boldsymbol{\S}$\textbf{11}.\quad  NUMERICAL EXAMPLES}
\setlength\parindent{2em}
\setcounter{theoremn}{0}
\renewcommand{\thepage}{\S11-\arabic{page}}
%%----------------------------------------------------------------------------------------------01

\ 
\\[-1.25cm]
\indent 

The basis for the calculations infra is $\S9$, $\#8$, which will not be quoted over and over, as well as the generalities in $\S10$, which will also be taken without attribution.
\\

\begin{x}{\small\bf EXAMPLE} \ %01
$\pi$ is $\ell$EL-computable, i.e., $\pi \in \ \R_{\ell\text{EL}}$.  
Thus write
\[
\frac{\pi}{4} 
\ = \ 
\sum\limits_{n=0}^\infty \hsx \frac{(-1)^n}{2n + 1}.
\]
Then the function 
\[
n \ra \frac{(-1)^n}{2n + 1}
\]
is $\ell$EL-computable.  As for convergence, the series
\[
\sum\limits_{n=0}^\infty \hsx \frac{(-1)^n}{2n + 1}
\]
is alternating, so $\forall \ x \in \ \N$ $(\xi(x) = x)$,
\begin{align*}
\abs{\sum\limits_{n= x + 1}^\infty \hsx  \frac{(-1)^n}{2n + 1}} \ 
&\leq \ 
\frac{1}{2(x + 1) + 1}
\\[11pt]
&=\ 
\frac{1}{2x + 3}
\\[11pt]
&\leq\ 
\frac{1}{x + 1}.
\end{align*}
These considerations establish the $\ell$EL-computability of $\ds\frac{\pi}{4}$, thus that of $\ds\pi = 4\big(\frac{\pi}{4}\big)$.
\end{x}

\begin{x}{\small\bf EXAMPLE} \ %02
$\elln(N)$ $(N = 1, 2, \ldots)$ is $\ell$EL-computable, i.e., $\elln(N)\in \ \R_{\ell\text{EL}}$.  
Thus write
\[
\elln\bigg(1 + \frac{1}{N}\bigg) 
\ = \ 
\sum\limits_{n=0}^\infty \hsx \frac{(-1)^n}{(n + 1)N^{n + 1}}
\]
%%----------------------------------------------------------------------------------------------02

\qquad\qquad $\implies$
\[
\elln(N + 1) 
\ = \ 
\elln(N) \ + \ \sum\limits_{n=0}^\infty \hsx \frac{(-1)^n}{(n + 1)N^{n + 1}}.
\]
Proceed by induction on \mN.  
When $N = 1$, $\elln(1) = 0$, which is obviously $\ell$EL-computable.  
So take $N > 1$ and suppose that $\elln(N)$ is $\ell$EL-computable. 
Since $R_{\ell\text{EL}}$ is a field, it need only be shown that
\[
\sum\limits_{n=0}^\infty \hsx \frac{(-1)^n}{(n + 1)N^{n + 1}} \in \ \R_{\ell\text{EL}}.
\]
But the function 
\[
n \ra \frac{(-1)^n}{(n + 1)N^{n + 1}}
\]
is $\ell$EL-computable.  
And $\forall \ x \in \ \N$ $(\xi(x) = x)$, 
\begin{align*}
\abs{\sum\limits_{n=x + 1}^\infty \hsx \frac{(-1)^n}{(n + 1)N^{n + 1}}} \ 
&\leq \ 
\frac{1}{(x + 2)N^{x + 2}}
\\[11pt]
&\leq \ 
\frac{1}{x + 1}.
\end{align*}
\end{x}

\begin{x}{\small\bf EXAMPLE} \ %03
Catalan's constant
\[
G 
\ = \ 
\sum\limits_{n=0}^\infty \hsx \frac{(-1)^n}{(2n + 1)^2}
\]
is $\ell$EL-computable, i.e., $G \in \ \R_{\ell\text{EL}}$.
\end{x}
\vspace{0.2cm}

\begin{x}{\small\bf EXAMPLE} \ %04
Euler's constant $\gamma$ is $\ell$EL-computable, i.e., $\gamma \in \ \R_{\ell\text{EL}}$.  
Thus write
\[
\gamma 
\ = \ 
\sum\limits_{n=0}^\infty \hsx \bigg(\frac{1}{n + 1} - \ell \bigg(1 + \frac{1}{n + 1}\bigg)\bigg).
\]
%%----------------------------------------------------------------------------------------------03
Then
\begin{align*}
\frac{1}{n + 1} - \elln\bigg(1 + \frac{1}{n + 1}\bigg) \ 
&=\ 
\sum\limits_{m=0}^\infty \hsx \frac{(-1)^m}{(m + 2)(n + 1)^{m + 2}}.
\end{align*}
Per $m$, 
\[
\alpha_n(m) 
\ = \ 
\frac{(-1)^m}{(m + 2)(n + 1)^{m + 2}}
\]
is a product of three bounded $\ell$EL-computable functions, thus $\alpha_n(m)$ is $\ell$EL-computable.  
As for convergence, $\forall \ m \in \ \N$ $(\xi(x) = x)$,
\begin{align*}
\abs{\sum\limits_{m = x + 1}^\infty \hsx \alpha_n(m)} \ 
&\leq \ 
\frac{1}{(x + 3)(n + 1)^{x + 3}}
\\[11pt]
&\leq \ 
\frac{1}{x + 3}
\\[11pt]
&\leq \  
\frac{1}{x + 1}.
\end{align*}
Therefore per $n$ the sum
\[
\Xi_n
\ = \ 
\sum\limits_{m = 0}^\infty \hsx \alpha_n(m)
\]
is $\ell$EL-computable. 
And finally
\begin{align*}
\abs{\sum\limits_{n = x + 1}^\infty \hsx \Xi_n} \ 
&\leq \ 
\sum\limits_{n = x + 1}^\infty \hsx \abs{\Xi_n}
\\[11pt]
&\leq \ 
\sum\limits_{n = x + 1}^\infty \  \abs{\sum\limits_{m = 0}^\infty \hsx \frac{(-1)^m}{(m + 2) (n + 1)^{m + 2}}}
\\[11pt]
&\leq \ 
\sum\limits_{n = x + 1}^\infty \hsx \frac{1}{2(n + 1)^2}
\\[11pt]
&= \ 
\frac{1}{2} \hsx \sum\limits_{n = x + 1}^\infty \hsx \frac{1}{(n + 1)^2}
%%----------------------------------------------------------------------------------------------04
\\[11pt]
&\leq \ 
\frac{1}{2} \hsx \sum\limits_{n = x + 1}^\infty \hsx \frac{1}{n(n + 1)}
\\[11pt]
&= \ 
\frac{1}{2} \hsx \sum\limits_{n = x + 1}^\infty \hsx \bigg(\frac{1}{n} - \frac{1}{n + 1}\bigg)
\\[11pt]
&\leq \ 
\frac{1}{2} \hsx \frac{1}{x + 1}
\\[11pt]
&\leq \ 
\frac{1}{x + 1}.
\end{align*}
Therefore
\[
\Xi 
\ = \ 
\sum\limits_{n = 0}^\infty \hsx \Xi_n 
\ = \ 
\gamma
\]
is $\ell$EL-computable.
\end{x}
\vspace{0.2cm}

\begin{x}{\small\bf EXAMPLE} \ %05
Liouville's number 
\[
L 
\ = \ 
\sum\limits_{n = 1}^\infty \hsx \frac{1}{10^{n!}}
\]
is $\ell$EL-computable, i.e., $L \in \ \R_{\ell\text{EL}}$.
\\[-0.25cm]

[As regards convergence, write
\[
L 
\ = \ 
\sum\limits_{n = 0}^\infty \hsx \frac{1}{10^{(n + 1)!}}
\]
and note that $\forall \ n \in \ \N$  $(\xi(x) = x)$,
\begin{align*}
\sum\limits_{n = x + 1}^\infty \hsx \frac{1}{10^{(n + 1)!}}\ 
&\leq 
\sum\limits_{n = x + 1}^\infty \hsx \frac{1}{(n + 1)!}
\\[11pt]
&= 
\frac{1}{(x + 2)!} + \frac{1}{(x + 3)!} +\frac{1}{(x + 4)!}  + \cdots
\\[11pt]
&= 
\frac{1}{(x + 2)!} \bigg(1 + \frac{1}{x + 3} + \frac{1}{(x + 3)(x + 4)} + \cdots \bigg)
\\[11pt]
&\leq 
\frac{1}{(x + 2)!} \bigg(1 + \frac{1}{2} + \frac{1}{2^2} + \cdots \bigg)
\\[11pt]
&= 
\frac{2}{(x + 2)!}
\\[11pt]
%%----------------------------------------------------------------------------------------------05
&\leq 
\frac{2}{2^{x + 1}}
\\[11pt]
&\leq 
\frac{1}{x + 1}.]
\\[11pt]
\end{align*}
\end{x}
\vspace{0.3cm}

\begin{x}{\small\bf EXAMPLE} \ %06
Let
\[
\zeta(x)
\ = \ 
\sum\limits_{n = 1}^\infty \hsx \frac{1}{n^x} \qquad (x > 1)
\]
and define
\[
f_\R:\N \ra \R
\]
by
\[
f_\R(k) \ = \ \zeta(k + 2) \qquad (k = 0, 1, \ldots).
\]
Then $f_\R$ is an $\ell$EL-computable function.  Consequently $\forall \ k$, 
\[
\zeta(k + 2) \ \in \ \R_{\ell\text{EL}}.
\]
In particular: $\zeta(3)$ is an $\ell$EL-computable real number.
\\[-0.25cm]

[Note: \ 
Put $\ell = k + 2$, thus $\ell = 2, 3, \ldots$ $-$then
\[
\zeta(\ell) 
\ = \ 
\int\limits_{1 > t_1 > \cdots > t_\ell > 0} \frac{\td t_1}{t_1} \cdots \frac{\td t_{\ell - 1}}{t_{\ell - 1}} \cdot \frac{\td t_\ell}{t_\ell}
\]
is a period, so
\[
\zeta(\ell) \ \in \ \R_{\ell\text{EL}} \qquad \text{(cf. $\S12$, $\#9$)},
\]
which is another way of looking at matters.]
\end{x}
\vspace{0.3cm}

\begin{x}{\small\bf FACTS} \ %07
\\[-0.25cm]

\hspace{0.65cm} \textbullet \quad $\forall \ k$, 
\[
f_\R(k + 1) \ < \ f_\R(k).
\]

\begin{align*}
[f_\R(k) -f_\R(k + 1)\ 
&=\ 
\sum\limits_{n=0}^\infty \hsx \frac{1}{(n + 1)^{k + 2}} - \sum\limits_{n=0}^\infty \hsx \frac{1}{(n + 1)^{k + 3}} 
\\[11pt]
%%----------------------------------------------------------------------------------------------06
&=\ 
\sum\limits_{n=0}^\infty \hsx\bigg(\frac{1}{(n + 1)^{k + 2}} -\frac{1}{(n + 1)^{k + 3}}\bigg)
\\[11pt]
&=\ 
\sum\limits_{n=0}^\infty \hsx \frac{1}{(n + 1)^{k + 2}}\bigg(1 - \frac{1}{n + 1}\bigg)
\\[11pt]
&=\ 
\sum\limits_{n=0}^\infty \hsx \frac{n}{(n + 1)^{k + 3}}
\\[11pt]
&>\ 
0.]
\end{align*}

\hspace{0.65cm} \textbullet \quad $\forall \ k$, 
\[
1 \ < \ f_\R(k) \ < \ 2.
\]

[$\ds f_\R(0) = \frac{\pi^2}{6} < 2$.  On the other hand, 
\[
f_\R(k) \ = \ 1 + \sum\limits_{n=2}^\infty \hsx \frac{1}{n^{k+2}} \ > \ 1.
\]
Therefore
\[
2
 \ > \ 
f_\R(0)
\ \geq \ 
f_\R(k)
 \ > \ 
 1.]
\]

\end{x}
\vspace{0.2cm}

\begin{x}{\small\bf EXAMPLE} \ %08
$\elln(\pi)$ is $\ell$EL-computable, i.e., $\elln(\pi) \in \ \R_{\ell\text{EL}}$.  
Thus write
\begin{align*}
\elln(\pi) - \elln(2) \ 
&=\ 
\sum\limits_{n=1}^\infty \hsx \frac{\zeta(2n)}{(2n)2^{2n - 1}}
\\[11pt]
&=\ 
\sum\limits_{n=0}^\infty \hsx \frac{\zeta(2(n + 1))}{2(n + 1)2^{2(n + 1) - 1}}
\\[11pt]
&=\ 
\sum\limits_{n=0}^\infty \hsx \frac{f_\R(2n)}{2(n + 1)2^{2n + 1}}.
\end{align*}
Then
\[
\elln(\pi) 
\ = \ 
 \elln(2) + \sum\limits_{n=0}^\infty \hsx \frac{f_\R(2n)}{2(n + 1)2^{2n + 1}}.
\]

%%----------------------------------------------------------------------------------------------07
Since $\elln(2) \in \ \R_{\ell\text{EL}}$ (cf. $\#2$), it suffices to examine the series
\[
\sum\limits_{n=0}^\infty \hsx \frac{f_\R(2n)}{2(n + 1) 2^{2n + 1}}.
\]
But

\hspace{0.5cm} (i) \quad 
$f_\R(n)$ is a bounded $\ell$EL-computable function, hence so is $f_\R(2n)$.
\\

\hspace{0.5cm} (ii) \quad
$\ds\frac{1}{2(n+1)} = \frac{1}{2n + 1 + 1}$ is an $\ell$EL-computable function.
\\

\hspace{0.5cm} (iii) \quad
$\ds\frac{1}{2^{2n + 1}}$ is an $\ell$EL-computable function.
\\[0.5cm]
Conclusion:
\[
\frac{f_\R(2n)}{2(n + 1) 2^{2n + 1}}
\]
is an $\ell$EL-computable function.

To handle the convergence,  $\forall \ x \in \hsx \N$  $(\xi(x) = x)$,
\begin{align*}
\sum\limits_{n=x + 1}^\infty \hsx \frac{f_\R(2n)}{2(n + 1) 2^{2n + 1}}\ 
&\leq \ 
\sum\limits_{n=x + 1}^\infty \hsx \frac{2}{2(n + 1) 2^{2n + 1}}
\\[11pt]
&\leq \ 
\sum\limits_{n=x + 1}^\infty \hsx \frac{1}{2^{2n + 1}}
\\[11pt]
&= \ 
\frac{1}{2^{2x + 3}} \cdot \frac{1}{1 - \frac{1}{4}}
\\[11pt]
&= \ 
\frac{4}{3 \cdot 2^{2x + 3}} 
\\[11pt]
&= \ 
\frac{1}{3 \cdot 2^{2x + 1}} 
%%----------------------------------------------------------------------------------------------08
\\[11pt]
&\leq \ 
\frac{1}{3(2x + 1)}
\\[11pt]
&= \ 
\frac{1}{6x + 3}
\\[11pt]
&\leq \ 
\frac{1}{x + 1}.
\end{align*}
\end{x}
\vspace{0.3cm}

\[
\text{APPENDIX}
\]
\vspace{0.3cm}

{\small\bf FACT} \quad 
If $\xi \in \ \R_{\ell\text{EL}}$, then $e^\xi \in \ \R_{\ell\text{EL}}$ 
and if $\xi \in \ \R_{\ell\text{EL}}$ is $> 0$, then $\elln(\xi) \in \ \R_{\ell\text{EL}}$. 
\\

{\small\bf FACT} \quad
If $\xi \in \ \R_{\ell\text{EL}}$, then
\[
\begin{cases}
\ \text{sin $\xi$} \\[5pt]
\ \text{cos $\xi$}
\end{cases}
\in \ \R_{\ell\text{EL}}.
\]
\\

{\small\bf FACT} \quad 
If $\xi \in \ \R_{\ell\text{EL}}$ and if $\abs{\xi} \leq 1$, then
\[
\begin{cases}
\ \text{Arc sin $\xi$} \\[5pt]
\ \text{Arc cos $\xi$}
\end{cases}
\in \ \R_{\ell\text{EL}}.
\]
\\

%%{\small\bf FACT} \quad 
%%If $\xi \in \ \R_{\ell\text{EL}}$, then
%%\[
%%\text{Arc tan $\xi \in \ \R_{\ell\text{EL}}$}.
%%\]

%%%%%%%%%%%%%%%%%%%%%%%%%%%%%%%%%%%%%%
%%%%%%%%%%%%%%%%%%%%%%%%%%%%%%%%%%%%%%
%%%%%%%%%%%%%%%%%%%%%%%%%%%%%%%%%%%%%%

%123
\chapter{
$\boldsymbol{\S}$\textbf{12}.\quad  $\R_{\sE^3}$ VERSUS $\R_{\sE^2}$}
\setlength\parindent{2em}
\setcounter{theoremn}{0}
\renewcommand{\thepage}{\S12-\arabic{page}}
%%----------------------------------------------------------------------------------------------01

\ 
\\[-1.25cm]
\indent 

Before getting down to business, there is a preliminary fact, frankly technical, which will be needed below.
\\

\begin{x}{\small\bf RAPPEL} \ %01
(cf. $\S 9$, $\# 5$) \ If $\alpha \in \R$ is $\sF$-computable, then there exist functions 
$f:\N \ra \N$ and $g:\N \ra \N$ in $\sF$ such that $\forall \ x$, 
\[
\abs{\frac{f(x) - g(x)}{x + 1} - \alpha} \ \leq \ \frac{1}{x + 1}.
\]

[Note: \ 
Here it is assumed that $\sF$ satisfies the standard conditions.]
\end{x}

\begin{x}{\small\bf SUBLEMMA} \ %02
If $\alpha \in \R$ is $\sF$-computable, then there exist functions $f_0$, $g_0$, $h_0$ in $\sF$ such that $\forall \ x$,
\[
\abs{\frac{f_0(x) - g_0(x)}{h_0(x) + 1} - \alpha} \ < \ \frac{1}{2(x + 1)}.
\]

PROOF \ 
Start with $u$, $v$, $w$ in $\sF$ such that $\forall \ x$,
\[
\abs{\frac{u(x) - v(x)}{w(x) + 1} - \alpha} \ \leq \ \frac{1}{x + 1}.
\]
Introduce
\[
\begin{cases}
\ f_0(x) \ = \ u(2(x + 1) + 1)\\[8pt]
\ g_0(x) \ = \ v(2(x + 1) + 1)\\[8pt]
\ h_0(x) \ = \ w(2(x + 1) + 1)
\end{cases}
.
\]
%%----------------------------------------------------------------------------------------------02
Then
\begin{align*}
\abs{\frac{f_0(x) - g_0(x)}{h_0(x) + 1} - \alpha} \ 
&=\ 
\abs{\frac{u(2(x+1) + 1) - v(2(x+ 1) + 1)}{h_0(2(x + 1) + 1)} - \alpha}
\\[11pt]
&\leq\ 
\frac{1}{2 (x + 1) + 1 + 1}
\\[11pt]
&=\ 
\frac{1}{2 x + 4)}
\\[11pt]
&=\ 
\frac{1}{2 (x + 2)}
\\[11pt]
&<\ 
\frac{1}{2 (x + 1)}.
\end{align*}
\end{x}

\begin{x}{\small\bf LEMMA} \ %03
 If $\alpha \in \R$ is $\sF$-computable, then there exist functions 
$f:\N \ra \N$ and $g:\N \ra \N$ in $\sF$ such that $\forall \ x$, 
\[
\abs{\frac{f(x) - g(x)}{x+ 1} - \alpha} 
\ < \ 
\frac{1}{x + 1}.
\]

PROOF \ 
Proceed as in $\S 9$, $\# 5$, taking $f, \ g \in \sF$ per $f_0$, $g_0$, $h_0$ to arrive at 
\begin{align*}
\abs{\frac{f(x) - g(x)}{x+ 1} - \alpha} \ 
&\leq \ 
\abs{\frac{f(x) - g(x)}{x+ 1} - \frac{f_0(x) - g_0(x)}{h_0(x) + 1}} + \abs{\frac{f_0(x) - g_0(x)}{h_0(x) + 1} - \alpha} 
\\[11pt]
&\leq 
\frac{1}{2(x + 1)} + \abs{\frac{f_0(x) - g_0(x)}{h_0(x) + 1} - \alpha} 
\\[11pt]
%%----------------------------------------------------------------------------------------------03
&< \ 
\frac{1}{2(x + 1)} + \frac{1}{2(x + 1)}
\\[11pt]
&= \ 
\frac{1}{x + 1}.
\end{align*}
\end{x}
\vspace{0.2cm}

\begin{x}{\small\bf \un{N.B.}} \ %04
The point, of course, is that in $\S 9$, $\# 5$, $\leq$ can be replaced by $<$.
\end{x}
\vspace{0.2cm}

\begin{x}{\small\bf LEMMA} \ %05
There exists a two-ary function $h$ in $\sE^3$ which is universal for the one-ary functions in $\sE^2$.
\\[-0.5cm]

[Note: \ 
I.e., the functions
\[
h \ \ra \ h(X,n) \qquad (X = 0, 1, 2, \ldots)
\]
exhaust the one-ary functions in $\sE^2$.]
\end{x}
\vspace{0.3cm}

\begin{x}{\small\bf RAPPEL} \ %06
The natural number function
\[
\text{pow}(x,y) \ = \ x^y
\]
is elementary (cf. $\S 6$, $\# 8$).
\end{x}
\vspace{0.3cm}

On general grounds, 
\[
\R_{\sE^2}
\ \subset \ 
\R_{\sE^3}
\]
and the claim is that this containment is strict.
\\[-0.5cm]

In detail: Start by defining a one-ary function $g$ as follows: 
$g(0) = 0$ and $\forall \ k  \in  \N$,
\[
g(k + 1) \ = \ 
\begin{cases}
\ 3 g(k) \quad \text{if} \quad 6 g(k) + 3 \leq h(k, 2 \cdot 3^{k + 1} - 1) \\[8pt]
\ 3 g(k) + 2 \qquad \text{otherwise}
\end{cases}
,
\]
thus
\[
g(k + 1) - 3 g(k) \in \ \{0, 2\}.
\]
\\

\begin{x}{\small\bf LEMMA} \ %07
$g \in \ \sE^3$.
\\[-0.25cm]
%%----------------------------------------------------------------------------------------------04

[Use the inequality $g(k) \leq 3^k - 1$.]
\end{x}
\vspace{0.2cm}

Put 
\[
\alpha 
\ = \ 
\sum\limits_{k = 0}^\infty \hsx \frac{g(k + 1)- 3 g(k)}{3^{k + 1}}.
\]
Then for any $K \in \ \N$, 
\[
\sum\limits_{k = 0}^{K - 1} \hsx \frac{g(k + 1)- 3 g(k)}{3^{k + 1}}
\ = \ 
\frac{g(K)}{3^K}
\]
$\implies$
\begin{align*}
0 \ 
&\leq \
\alpha - \frac{g(K)}{3^K}
\\[11pt]
&= \
\sum\limits_{k = K}^\infty \hsx \frac{g(k + 1)- 3 g(k)}{3^{k + 1}}
\\[11pt]
&\leq \
2 \sum\limits_{k = K}^\infty \hsx \frac{1}{3^{k + 1}}
\\[11pt]
&= \
2 \bigg(\frac{1}{3^{K + 1}} + \frac{1}{3^{K + 2}}  + \frac{1}{3^{K + 3}}  + \cdots \bigg)
\\[11pt]
&= \
\frac{2}{3^{K + 1}} \bigg(1 + \frac{1}{3}  + \frac{1}{3^2}  + \cdots \bigg)
\\[11pt]
&= \
\frac{2}{3^{K + 1}} \cdot \frac{3}{2}
\\[11pt]
&= \
\frac{1}{3^K}
\end{align*}
\qquad \qquad $\implies$
\begin{align*}
0 \ 
&\leq \
\alpha - \frac{g(K)}{3^K}
\\[11pt]
&= \
\alpha - \frac{g(K)}{(3^K - 1) + 1}
\\[11pt]
&\leq \
\frac{1}{3^K}
\\[11pt]
&\leq \
\frac{1}{K + 1}.
\end{align*}
%%----------------------------------------------------------------------------------------------05
Therefore the real number $\alpha$ is $\sE^3$-computable:
\[
\alpha \in \ \R_{\sE^3}.
\]
Still
\[
\alpha \not\in \ \R_{\sE^2}.
\]

To establish this, proceed by contradiction and assume that $\forall \ x \in  \N$, 
\[
\abs{\frac{f(x) - g(x)}{x + 1} - \alpha}
\ < \ 
\frac{1}{x + 1} \qquad \text{(cf. $\# 3$).}
\]
Here $f$ and $g$ belong to $\sE^2$, as does $\abs{f - g}$.  
And, since $\alpha$ is nonnegative, $\forall \ x \in  \N$,
\[
\abs{\frac{\abs{f(x) - g(x)}}{x + 1} - \alpha}
\ < \ 
\frac{1}{x + 1}.
\]

PROOF \ 
For all real \mS and \mT,
\[
\abs{\abs{S} - \abs{T}} \ \leq \ \abs{S - T}
\]
\qquad \qquad $\implies$
\begin{align*}
\abs{\frac{\abs{f(x) - g(x)}}{\abs{x + 1}} - \abs{\alpha}} \ 
&=\ 
\abs{\frac{\abs{f(x) - g(x)}}{x + 1} - \alpha}
\\[11pt]
&\leq\ 
\abs{\frac{f(x) - g(x)}{x + 1} - \alpha}
\\[11pt]
&< \ 
\frac{1}{x + 1}.
\end{align*}
Choose now per $\# 5$ a natural number \mX such that $\forall \ x \in \N$, 
\[
\abs{f(x) - g(x)} \ = \ h(X,x),
\]
hence
\[
\abs{\frac{h(X,x)}{x + 1} - \alpha} \ < \ \frac{1}{x + 1}.
\]
%%----------------------------------------------------------------------------------------------06
In particular: 
\[
\abs{\frac{h(X,2 \cdot 3^{X + 1} - 1)}{2 \cdot 3^{X + 1}} - \alpha} 
\ < \ 
\frac{1}{2 \cdot 3^{X + 1}}.
\]
There are then two possibilities, both of which lead to a contradiction.

\hspace{0.5cm} \textbullet \quad
\begin{align*}
&6 g(X) + 3 \ \leq \ h(X,2 \cdot 3^{X + 1} - 1)
\\[11pt]
\implies \hspace{1cm} 
&
\\[11pt]
&\frac{6 g(X) + 3}{2 \cdot 3^{X + 1}} \leq \frac{h(X,2 \cdot 3^{X + 1} - 1)}{2 \cdot 3^{X + 1}}
\\[11pt]
\implies \hspace{1cm} 
&
\\[11pt]
&\frac{3 g(X)}{3^{X + 1}} + \frac{1}{2 \cdot 3^X}  \ \leq \ \frac{h(X,2 \cdot 3^{X + 1} - 1)}{2 \cdot 3^{X + 1}}
\\[11pt]
\implies \hspace{1cm} 
&
\\[11pt]
&\frac{g(X)}{3^X} + \frac{1}{2 \cdot 3^X}  \ \leq \ \frac{h(X,2 \cdot 3^{X + 1} - 1)}{2 \cdot 3^{X + 1}}
\\[11pt]
&
\hspace{1cm} < \ \alpha + \frac{1}{2 \cdot 3^{X + 1}} .
\end{align*}
But
\begin{align*}
0 \ 
&\leq \ 
\alpha - \frac{g(X + 1)}{3^{X + 1}} 
\\[11pt]
&\leq\ 
\frac{1}{3^{X + 1}}
\\[11pt]
\implies \hspace{1cm} 
&
\\[11pt]
\alpha \ 
&\leq \ 
\frac{g(X + 1)}{3^{X + 1}} + \frac{1}{3^{X + 1}}
\\[11pt]
&= \ 
\frac{3 g(X)}{3^{X + 1}}
\\[11pt]
&= \ 
\frac{g(X)}{3^X} + \frac{1}{3^{X + 1}}
\\[11pt]
%%----------------------------------------------------------------------------------------------07
\implies \hspace{1cm} 
&
\\[11pt]
\cdots \ 
&< \ 
\alpha + \frac{1}{2 \cdot 3^{X + 1}}
\\[11pt]
&\leq \ 
\frac{g(X)}{3^X} + \frac{1}{3^{X + 1}} + \frac{1}{2 \cdot 3^{X + 1}}
\\[11pt]
&= \ 
\frac{g(X)}{3^X} + \frac{1}{3^{X + 1}}  \bigg(1 + \frac{1}{2}\bigg)
\\[11pt]
&= \ 
\frac{g(X)}{3^X} + \frac{1}{3^{X + 1}}  \bigg(\frac{3}{2}\bigg)
\\[11pt]
&= \ 
\frac{g(X)}{3^X} + \frac{1}{2 \cdot 3^{X}}.
\\[11pt]
\end{align*}
Contradiction.
\\

\hspace{0.5cm} \textbullet \quad
\begin{align*}
&6 g(X) + 3 \ 
> \ 
h(X,2 \cdot 3^{X + 1} - 1)
\\[11pt]
\implies \hspace{1cm} 
&
\\[11pt]
&\frac{g(X)}{3^X} + \frac{1}{2 \cdot 3^X} \ > \ \frac{h(X, 2 \cdot 3^{X + 1} - 1)}{2 \cdot 3^{X + 1}}
\\[11pt]
&\hspace{2.6cm} > \ \alpha - \frac{1}{2 \cdot 3^{X + 1}}.
\end{align*}
But
\begin{align*}
\alpha \ 
&\geq \ 
\frac{g(X + 1)}{3^{X + 1}}
\\[11pt]
\implies \hspace{1cm} 
&
\\[11pt]
\alpha \ 
&=\ 
\frac{3 g(X) + 2}{3^{X + 1}}
\\[11pt]
&=\ 
\frac{g(X)}{3^X} + \frac{2}{3^{X + 1}}
\\[11pt]
\implies \hspace{1cm} 
&
\\[11pt]
\cdots \ 
&> \ 
\alpha - \frac{1}{2 \cdot 3^{X + 1}}
%%----------------------------------------------------------------------------------------------08
\\[11pt]
&\geq \ 
\frac{g(X)}{3^X} + \frac{2}{3^{X + 1}} - \frac{1}{2 \cdot 3^{X + 1}}
\\[11pt]
&= \ 
\frac{g(X)}{3^X} + \frac{1}{3^{X + 1}} \bigg(2 - \frac{1}{2}\bigg)
\\[11pt]
&= \ 
\frac{g(X)}{3^X} + \frac{1}{3^{X + 1}} \bigg(\frac{3}{2}\bigg)
\\[11pt]
&= \ 
\frac{g(X)}{3^X} + \frac{1}{2 \cdot 3^{X}}.
\end{align*}
Contradiction.

\begin{x}{\small\bf REMARK} \ %08
Suppose that $n \geq 2$ $-$then it can be shown that $\# 5$ remains valid if $\sE^3$ is replaced by $\sE^{n+1}$ and $\sE^2$ is replaced by $\sE^n$.  
This said, the argument above goes through without change, the conclusion being that 
\[
\exists \ \alpha \in \ \R_{\sE^{n + 1}} : \alpha \not\in \R_{\sE^{n}}.
\]
\end{x}
\vspace{0.3cm}

\begin{x}{\small\bf THEOREM} \ %09
\[
\sP_\text{KZ} \ \subset \ \R_\text{$\ell$EL}.
\]

[For the details, see Katvin Tent and Martin Zeigler
\footnote[2]{\textit{M\"unster Journal of Mathematics} {\bf{3}} (2010), pp. 43-66.}
.]

So there is a chain
\[
\sP_\text{KZ} 
\ \subset \ 
\R_\text{$\ell$EL}
\ \subset \ 
\R_{\sE^2}
\ \subset \ 
\R_{\sE^3}
\ \equiv \ 
\R_\text{EL}.
\]
And in view of what has been said above, the containment
\[
\R_{\sE^2}
\ \subset \ 
\R_{\sE^3}
\]
is strict.
\end{x}
\vspace{0.3cm}

%%----------------------------------------------------------------------------------------------09
\begin{x}{\small\bf \un{N.B.}}  \ %10
It is unknown whether $\ell$EL equals $\sE^2$ or not.
\end{x}
\vspace{0.3cm}

\begin{x}{\small\bf EXAMPLE} \ %11
$\pi \in \ \R_{\text{$\ell$EL}}$ (cf. $\S 11$, $\# 1$) but actually $\pi \in \ \sP_{\text{KZ}}$ (cf. $\S 2$, $\# 6$).
\\[-0.25cm]

[Note: \ 
$e \in \ \R_{\text{$\ell$EL}}$ but it is not known if  $e \in \ \sP_{\text{KZ}}.]$ 
\end{x}
\vspace{0.3cm}

\begin{x}{\small\bf EXAMPLE} \ %12
$e^\pi \in \ R_{\text{$\ell$EL}}$ (see the Appendix to $\S 11$).
\end{x}
\vspace{0.3cm}

%%%%%%%%%%%%%%%%%%%%%%%%%%%%%%%%%%%%%%
%%%%%%%%%%%%%%%%%%%%%%%%%%%%%%%%%%%%%%
%%%%%%%%%%%%%%%%%%%%%%%%%%%%%%%%%%%%%%

%12
\chapter{
$\boldsymbol{\S}$\textbf{13}.\quad  RECURSIVE FUNCTIONS}
\setlength\parindent{2em}
\setcounter{theoremn}{0}
\renewcommand{\thepage}{\S13-\arabic{page}}
%%----------------------------------------------------------------------------------------------01

\begin{x}{\small\bf DEFINITION} \ %01
The set of \un{recursive functions}, denoted by $\fR$, is the inductive closure of the initial functions with respect to the operations of composition, 
primitive recursion, and minimization.
\end{x}

\begin{x}{\small\bf \un{N.B.}} \ %02
Obviously
\[
\fP\fR \ \subset \ \fR.
\]

[Note: \ 
The containment is proper (the Ackermann function figuring in the Appendix to $\S5$ is recursive but not primitive recursive).]
\end{x}

\begin{x}{\small\bf REMARK} \ %03
An important property of $\fP\fR$ is that it is a recursively enumerable subset of $\fR$, i.e., there is a two-way function $u(m,n)$ that enumerates 
the primitive recursive functions in the sense that
\\[-0.25cm]

\hspace{0.5cm} \textbullet \quad 
$\forall \ f \in \ \fP\fR$, $\exists \ m:$ $\forall \ n$, $f(n) = u(m,n)$.
\\[-0.25cm]

\hspace{0.5cm} \textbullet \quad 
$\forall \ m$, $u(m,-) \in \ \fP\fR$.
\\[-0.25cm]

[Note: \ 
On the other hand, $\fR$ itself is not recursively enumerable.]
\end{x}
\vspace{0.2cm}

In the theory developed in $\S8$, take $\sF = \fR$ (the standard conditions are then obviously in force).  
So an $\fR$-sequence is a function $A:\N \ra \Q$ that has a representation of the form
\[
A(x) 
\ = \ 
\frac{f(x) - g(x)}{h(x) + 1} \qquad (x = 0, 1, 2, \ldots),
\]
where $f$, $g$, $h:\N \ra \N$ belong to $\fR$ and a real number $\alpha$ is said to be \un{$\fR$-computable} 
if there exists an $\fR$-sequence $A:\N \ra \Q$ such that $\forall \ x$,
\[
\abs{A(x) - \alpha}
\ \leq \ 
\frac{1}{x + 1}.
\]

%%----------------------------------------------------------------------------------------------02

\begin{x}{\small\bf \un{N.B.}} \ %04
Rather than working with $(x + 1)^{-1}$one can work instead with $2^{-x}$, 
either definition leading to the same set of $\fR$-computable real numbers.
\\[-0.5cm]

[Note: \ 
This is not always permissible.  
E.g.: Take $\sF = \sE^2$ $-$then the use of $2^{-x}$ would imply that the $\sE^2$-computables are precisely the rationals, which is untenable.  
However, the switch to $2^{-x}$ is permissible if $\sF = \sE^n$ $(n \geq 3)$, in particular if $n = 3$ $(\implies \sE^3 = \text{EL})$ 
or if $\sF = \fP\fR$.]
\end{x}
\vspace{0.2cm}

\begin{x}{\small\bf NOTATION} \ %05
Denote the set of all $\fR$-computable real numbers by the symbol $\R_\fR$ (cf. $\S8$, $\#10$).
\end{x}
\vspace{0.3cm}

\begin{x}{\small\bf THEOREM} \ %06
$\R_\fR$ is a real closed field (cf. $\S8$, $\#21$).
\\[-0.5cm]

[Note: \ 
In addition, $\R_\fR$ is countable.]
\end{x}
\vspace{0.3cm}

\begin{x}{\small\bf \un{N.B.}} \ %07
It is customary to refer to the elements of $\R_\fR$ as simply the \un{computable reals}.
\end{x}
\vspace{0.2cm}

\begin{x}{\small\bf EXAMPLE} \ %08
\[
\sP_\text{KZ} \ \subset \ \R_{\ell\text{EL}} \ \text{(cf. $\S12$, $\#9$)} \ \subset \ \R_{\fP\fR} \ \subset \ \R_\fR.
\]
Therefore periods are computable.
\end{x}
\vspace{0.3cm}

\begin{x}{\small\bf EXAMPLE} \ %09
Chaitins constant(s) $\Omega$ is (are) not computable.
\end{x}
\vspace{0.3cm}

While the very definition of ``computable real'' involves recursive functions, matters can also be formulated in terms of primitive recursive functions.

\begin{x}{\small\bf DEFINITION} \ %10  
\quad
Let $\alpha$ be a real number \ $-$then \ a \ \un{primitive recursive} \un{approximation} 
of $\alpha$ is a pair $(A,E)$ of $\fP\fR$-sequences $A$, $E:\N \ra \Q$ such that \mE is
%%----------------------------------------------------------------------------------------------03
monotonically decreasing to 0 and such that $\forall \ x$, 
\[
\abs{A(x) - \alpha} 
\ \leq \ 
E(x).
\]

[Note: \ In general, \mE depends on $\alpha$.]
\end{x}
\vspace{0.3cm}

\begin{x}{\small\bf EXAMPLE} \ %11
Suppose that $\alpha \in \ \R_{\fP\fR}$, so there exists a $\fP\fR$-sequence $A:\N \ra \Q$ such that $\forall \ x$, 
\[
\abs{A(x) - \alpha} 
\ \leq \ 
\frac{1}{x + 1}.
\]
Then the pair
\[
\bigg(A(x),\frac{1}{x + 1}\bigg)
\]
is a primitive recursive approximation to $\alpha$.
\end{x}
\vspace{0.3cm}

\begin{x}{\small\bf THEOREM} \ %12
A real number $\alpha$ is computable iff it has a primitive recursive approximation.
\end{x}
\vspace{0.3cm}

One direction is straightforward.  
Thus consider a real number $\alpha$ with the stated property.  
Define $s:\N \ra \N$ by the rule
\[
s(x)
\ = \ 
\min\bigg\{n : E(n) \leq \frac{1}{x + 1}\bigg\}.
\]
Then $s$ is recursive, hence $A \circ s$ is an $\fR$-sequence (cf. $\S8$, $\#13$) and $\forall \ x$,
\[
\abs{A(s(x)) - \alpha} 
\ \leq \ 
\frac{1}{x + 1}.
\]
Therefore $\alpha$ is $\fR$-computable.  
\\[-0.25cm]

In the other direction:
\\[-0.25cm]

\begin{x}{\small\bf LEMMA} \ %13
Suppose that $\alpha \in \ \R_\fR$ $-$then there exists a pair $(A,E)$ of $\fP\fR$-sequences $A$, $E:\N \ra \Q$ with the property that there are elements of
$E(\N)$ which are arbitrarily close to 0 and such that $\forall \ x$,
\[
\abs{A(x) - \alpha} 
\ \leq \ 
E(x).
\]
%%----------------------------------------------------------------------------------------------04

PROOF \ 
The assumption on $\alpha$ implies that there exists an $\fR$-sequence $A^\prime$ such that $\forall \ x$,
\[
\abs{A^\prime(x) - \alpha} 
\ \leq \ 
\frac{1}{x + 1}.
\]
This said, choose a surjective primitive recursive function $f:\N \ra \N$ such that 
\[
A(x) 
\ = \ 
A^\prime(f(x))
\]
is a $\fP\fR$-sequence.  Put
\[
E(x) 
\ = \ 
\frac{1}{f(x) + 1}.
\]
Then
\begin{align*}
\abs{A(x) - \alpha}  \ 
&=\ 
\abs{A^\prime(f(x)) - \alpha} 
\\[11pt]
&\leq 
\frac{1}{f(x) + 1} 
\\[11pt]
&=\ 
E(x).
\end{align*}

To finish the proof of $\#12$, one has only to take the data supplied by $\#13$ and transform it into that required of $\#10$.  Using primes, put
\[
\begin{cases}
\ E^\prime(n) \ = \ \min\{E(i) : 0 \leq i \leq n\} \\[8pt]
\ k(n) \ = \ \min\{i : 0 \leq i \leq n, \ E(i) = E^\prime(n)\} \\[8pt]
\ A^\prime(n) \ = \ A(k(n))
\end{cases}
.
\]
Then the pair $(A^\prime, E^\prime)$ is a primitive recursive approximation of $\alpha$.
\end{x}
\vspace{0.3cm}

\begin{x}{\small\bf RAPPEL} \ %14
To say that a real number is computable means that there exists an $\fR$-sequence $A:\N \ra \Q$ such that $\forall \ x$, 
\[
\abs{A(x) - \alpha} 
\ \leq \ 
\frac{1}{x + 1}.
\]
\end{x}
\vspace{0.3cm}

Question: Can one instead utilize a $\fP\fR$-sequence?  The answer in general is ``no''.
\\
%%----------------------------------------------------------------------------------------------05

\begin{x}{\small\bf EXAMPLE} \ %15
Let $f:\N \ra \{0,1\}$ be recursive and put
\[
\alpha 
\ = \ 
\sum\limits_{n=0}^\infty \hsx \frac{f(n)}{4^n}.
\]
Then $\alpha \in \ \R_\fR$ and the claim is that there is no $\fP\fR$-sequence $A:\N \ra \Q$ with the property that $\forall \ x$, 
\[
\abs{A(x) - \alpha}
\ \leq \ 
\frac{1}{x + 1}.
\]

[The initial observation is that if $k \in \N$, $q \in \Q$, and 
\[
\abs{q - \alpha}
\ \leq \ 
\frac{1}{4^{k+1}},
\]
then 
\[
f(k) 
\ = \ 
\bigg[\frac{[2 \cdot 4^k q + 1/2] \hsx \modx 4}{2}\bigg] \qquad \text{(cf. infra).}
\]
Granted this, consider a $\fP\fR$-sequence $A:\N \ra \Q$ subject to $\forall \ x$,
\[
\abs{A(x) - \alpha}
\ \leq \ 
\frac{1}{x + 1},
\]
and let
\[
q 
\ = \ 
A(4^{k+1} - 1),
\]
thus
\[
\abs{q - \alpha}
\ \leq \ 
\frac{1}{4^{k+1} - 1 + 1}
\ = \ 
\frac{1}{4^{k+1}},
\]
and the formula for $f(k)$ implies that $f$ is primitive recursive.  
Accordingly, if $\alpha$ is constructed by using a function $f:\N \ra \{0,1\}$ that is recursive but not primitive recursive, i.e., if
\[
f \in \ \fR 
\quad \text{but} \quad 
f \not\in \fP\fR,
\]
then there will be no $\fP\fR$-sequence $A:\N \ra \Q$ per supra.]
%%----------------------------------------------------------------------------------------------06

[Details:
\[
2 \alpha 
\ = \ 
\sum\limits_{n=0}^\infty \hsx \frac{2 f(n)}{4^n}
\]
$\implies$
\[
[2 \cdot 4^k \alpha] \hsx \modx 4 
\ = \ 
2 f(k).
\]
And
\[
2 \cdot 4^k \alpha 
\ \leq \ 
2 \cdot 4^k q + 1/2 
\ \leq \ 
2 \cdot 4^k \alpha + 1
\]
$\implies$
\[
[2 \cdot 4^k q + 1/2 ] 
\ = \ 
[2 \cdot 4^k \alpha] + d \qquad (d = 0 \ \text{or} \ d = 1)
\]
$\implies$
\[
[2 \cdot 4^k q + 1/2] \hsx \modx 4 
\ = \ 
2 f(k) + d.]
\]

\end{x}
\vspace{0.3cm}

\begin{x}{\small\bf DEFINITION} \ %16 
\ 
Let $\alpha$ be a real number $-$then $\alpha$ has a 
\un{primitive recursive} \un{nested interval representation} if there are $\fP\fR$-sequences $f$, $g:\N \ra \Q$ such that $\forall \ x$, 
\[
f(x) 
\ \leq \ 
f(x + 1) 
\ \leq \ 
\alpha 
\ \leq \ 
g(x + 1) 
\ \leq \ 
g(x)
\]
and 
\[
\lim\limits_{x \ra \infty} (g(x) - f(x)) \ = \ 0.
\]
\end{x}
\vspace{0.3cm}

\begin{x}{\small\bf LEMMA} \ %17
A real number $\alpha$ has a primitive recursive nested interval representation iff it admits a primitive recursive approximation.
\\[-0.25cm]

PROOF \ 

\qquad $\implies$ \quad Given $f$, $g$, let
\[
A(x) 
\ = \ 
\frac{g(x) + f(x)}{2} 
\quad \text{and} \quad
E(x) 
\ = \ 
\frac{g(x) - f(x)}{2}.
\]
%%----------------------------------------------------------------------------------------------07
Then
\begin{align*}
\abs{A(x) - \alpha} \ 
&=\ 
\abs{\frac{g(x) + f(x)}{2} - \alpha} 
\\[11pt]
&= 
\abs{\frac{g(x)}{2} - \frac{\alpha}{2} + \frac{f(x)}{2} - \frac{\alpha}{2}}
\\[11pt]
&\leq 
\abs{\frac{g(x)}{2} - \frac{\alpha}{2}} + \abs{\frac{f(x)}{2} - \frac{\alpha}{2}}
\\[11pt]
&=
\frac{g(x)}{2} - \frac{\alpha}{2} + \frac{\alpha}{2} - \frac{f(x)}{2}
\\[11pt]
&= 
\frac{g(x) - f(x)}{2}
\\[11pt]
&= 
E(x).
\end{align*}
And
\[
E(x + 1) \ \leq \ E(x)
\]
iff
\[
g(x + 1) - f(x + 1)  
\ \leq \ 
g(x) - f(x)
\]
iff
\[
g(x + 1) - g(x)  
\ \leq \ 
f(x + 1) - f(x)
\]
iff
\[
g(x) - g(x + 1) 
\ \geq \ 
f(x) - f(x + 1).
\]
But
\[
\begin{cases}
\ g(x) - g(x + 1)  \geq 0 \\[11pt]
\ f(x) - f(x + 1) \leq 0
\end{cases}
. \hspace{2.7cm}
\]
\\[-0.5cm]

\qquad $\impliedby$ Given \mA, \mE, let
\[
\begin{cases}
\ f(x) \ = \ \max\{A(n) - E(n): n \leq x\} \\[11pt]
\ g(x) \ = \ \min\{A(n) + E(n) : n \leq x\}
\end{cases}
.
\]
\end{x}
\vspace{0.3cm}
%%----------------------------------------------------------------------------------------------08

\begin{x}{\small\bf SCHOLIUM} \ %18
A real number $\alpha$ is computable iff it has a primitive recursive nested interval representation.
\end{x}
\vspace{0.3cm}

\begin{x}{\small\bf EXAMPLE} \ %19
\[
\Q \ \subset \ \R_\fR.
\]

[Given $q \in \ \Q$, let
\[
\begin{cases}
\ f(x) = q - 2^{-x}\\[8pt]
\ g(x) = q + 2^x
\end{cases}
.]
\]
\end{x}
\vspace{0.3cm}

%%%%%%%%%%%%%%%%%%%%%%%%%%%%%%%%%%%%%%
%%%%%%%%%%%%%%%%%%%%%%%%%%%%%%%%%%%%%%
%%%%%%%%%%%%%%%%%%%%%%%%%%%%%%%%%%%%%%

%12
\chapter{
$\boldsymbol{\S}$\textbf{14}.\quad  EXPANSION THEORY}
\setlength\parindent{2em}
\setcounter{theoremn}{0}
\renewcommand{\thepage}{\S14-\arabic{page}}
%%----------------------------------------------------------------------------------------------01

\ 
\\[-1.25cm]
\indent 

Let $b > 1$ be a natural number and let $\alpha$ be a nonnegative real number $-$then 
$\alpha$ has a \un{$b$-adic representation} if there exists a recursive function 
$f:\N \ra \{0, 1, \ldots, b-1\}$ such that 
\[
\alpha \ = \  \sum\limits_{n = 0}^\infty \hsx \frac{f(n)}{b^n}.
\]

\begin{x}{\small\bf LEMMA} \ %01
If $\alpha$ has a $b$-adic representation, then $\alpha \in \R_{\fR}$.
\\[-.3cm]

PROOF \ 
Put
\[
A(x) 
\ = \ 
\sum\limits_{n = 0}^{x + 1} \hsx \frac{f(n)}{b^n} \qquad (x = 0, 1, 2, \ldots).
\]
Then $A:\N \ra \Q$ is an $\R$-sequence and 
\begin{align*}
\abs{A(x) - \alpha} \ 
&=\ 
\abs{\sum\limits_{n = 0}^{x + 1} \hsx \frac{f(n)}{b^n} - \sum\limits_{n = 0}^\infty \hsx \frac{f(n)}{b^n}}
\\[11pt]
&=\ 
\sum\limits_{n = x + 2}^\infty \hsx \frac{f(n)}{b^n}
\\[11pt]
&\leq\ 
\sum\limits_{n = x + 2}^\infty \hsx \frac{b - 1}{b^n}
\\[11pt]
&\leq\ 
b \hsx \sum\limits_{n = x + 2}^\infty \hsx \frac{1}{b^n}
\\[11pt]
&=\ 
b \bigg(\frac{1}{b^{x + 2}} + \frac{1}{b^{x + 3}} +\cdots \bigg)
\\[11pt]
&=\ 
\frac{b}{b^{x + 2}} \bigg(1 +  \frac{1}{b} + \cdots \bigg)
%%----------------------------------------------------------------------------------------------02
\\[11pt]
&=\ 
\frac{b}{b^{x + 2}} \frac{1}{1 - \frac{1}{b}}
\\[11pt]
&=\ 
\frac{b}{b^{x + 2}} \frac{b}{b - 1}
\\[11pt]
&=\ 
\frac{b^2}{b^{x + 2}} \frac{1}{b - 1}
\\[11pt]
&=\ 
\frac{1}{b^x} \frac{1}{b - 1}
\\[11pt]
&\leq\ 
\frac{1}{b^x} \qquad (2 \leq b \implies 1 \leq b - 1 \implies \frac{1}{b-1} \leq 1)
\\[11pt]
&=\ 
b^{-x}
\\[11pt]
&\leq\ 
2^{-x}.
\end{align*}
\end{x}

\begin{x}{\small\bf NOTATION} \ %02
$\R_{\fR}^b$ is the set of nonnegative real numbers which admit a $b$-adic representation.
\end{x}

Therefore
\[
\R_\fR^b \ \subset \ \R_\fR \hsx \cap \hsx [0,\infty[ \hsx.
\]

\begin{x}{\small\bf DEFINITION} \ %03
A subset $A \subset \ \N$ is \un{recursive} if its characteristic function $\chi_A$ is recursive, i.e., if $\chi_A \in \ \fR$.
\end{x}
\vspace{0.2cm}

\begin{x}{\small\bf EXAMPLE} \ %04
Suppose that $\alpha \in \R_{\fR}^2$, hence
\[
\alpha 
\ = \ \sum\limits_{n = 0}^\infty \ \frac{f(n)}{2^n} \qquad (f(n) \in \{0,1\}).
\]
Let
\[
A 
\ = \ 
\{n : f(n) = 1\}.
\]
%%----------------------------------------------------------------------------------------------03
Then
\[
\alpha 
\ = \ \sum\limits_{n \in A} 2^{-n}
\]
is computable and $\chi_A = f$, so \mA is recursive (this being the case of $f$).
\end{x}
\vspace{0.2cm}

\begin{x}{\small\bf SUBLEMMA} \ %05
Suppose that $\alpha$ is an irrational computable real number, thus there exists an $\fR$-sequence 
$A:\N \ra \Q$ such that $\forall \ x$, 
\[
\abs{A(x) - \alpha}
\ \leq \ 
2^{-x} \qquad \text{(cf. $\S13$, $\#4$).}
\]
Let $f:\N \ra \Q$ be an $\R$-sequence $-$then $\forall \ n$ $\exists \ x$ such that 
\[
\abs{A(x) - f(n)}
\ > \ 
2^{-x}.
\]
And
\[
\begin{cases}
\ A(x) - f(n) > 2^{-x} \implies x > f(n) \\[8pt]
\ A(x) - f(n) < 2^{-x} \implies x < f(n)
\end{cases}
.
\]
\end{x}
\vspace{0.3cm}

\begin{x}{\small\bf THEOREM} \ %06
\[
\R_\fR^b 
\ = \ 
\R_\fR \hsx \cap \hsx [0,\infty[\hsx .
\]

PROOF \ 
Take an $\alpha$ in $\R_\fR \hsx \cap \hsx [0,\infty[$.  
If it is rational, say $\ds\frac{p}{q}$, work relative to the base $b$ and carry out the long division of $p$ by $q$.  
Assume, therefore, that $\alpha$ is irrational and without loss of generality take $\alpha$ in $\R_\fR \  \cap \  ]0,1[$.  
Suppose now that we have found natural numbers $n_0 = 0, \ldots, n_k$ such that 
$0 \leq n_j \leq b -1$ $(1 \leq j \leq k)$ and such that
\[
\sum\limits_{j = 0}^k \hsx n_j b^{-j} 
\ < \ 
\alpha
\ < \ 
\sum\limits_{j = 0}^{k-1} \hsx n_j b^{-j} + (n_k + 1) b^{-k}.
\]
Applying $\#5$, there is a unique $x \in \  \{0, 1, \ldots, b-1\}$ for which
%%----------------------------------------------------------------------------------------------04
\[
\sum\limits_{j = 0}^k \hsx n_j b^{-j} + x b ^{-k - 1}
\ < \
\alpha
\ < \
\sum\limits_{j = 0}^k \hsx n_j b^{-j} + (x + 1) b^{-k - 1}.
\]
Definition: $n_{k+1} = x$.  With this agreement, set $f(k) = n_k$ $-$then $f:\N \ra \N$ is a recursive function, 
$f(k) \in \ \{0, 1, \ldots, b-1\}$, and
\[
\alpha 
\ = \ 
\sum\limits_{k = 0}^\infty \hsx \frac{f(k)}{b^k}.
\]
\end{x}
\vspace{0.3cm}

\begin{x}{\small\bf \un{N.B.}}  \ %07
It is to be emphasized that the realization
\[
\R_\fR^b \ = \ \R_\fR \hsx \cap \hsx [0,\infty[
\]
is valid for all $b$.
\end{x}
\vspace{0.2cm}

\begin{x}{\small\bf EXAMPLE} \ %08
Given a subset $A \subset \  \N$, put
\[
\alpha_A
\ = \ 
\sum\limits_{n \in \A} \ 2^{-n}.
\]
Then $\alpha_A$ is computable iff \mA is recursive.
\\[-0.25cm]

PROOF \ 
If $\alpha_A$ is computable, then $\alpha_A \in \ \R_\fR^2$ (cf. $\# 6$), so \mA is recursive (cf. $\# 4$).
   
\noindent On the other hand, if \mA is recursive, let
\[
s_x 
\ = \ 
\sum\limits_{n \leq x, \hsx n \in \ \A} \ 2^{-n}.
\]
Take
\[
\begin{cases}
\ f(x) = s_x \\[3pt]
\ g(x) = s_x + 2^{-x}
\end{cases}
.
\]
Then
\begin{align*}
s_x \ 
&\leq \ 
\alpha_A 
\\[11pt]
&=\ 
\sum\limits_{n \in \A} \ 2^{-n}
%%----------------------------------------------------------------------------------------------05
\\[11pt]
&=\ 
s_x + \sum\limits_{n > x, \hsx n \in \A} \ 2^{-n}
\\[11pt]
&\leq\ 
s_x + 2^{-(x + 1)} + 2^{-(x + 2)} + \cdots
\\[11pt]
&=\ 
s_x + 2^{-(x + 1)}(1 + 2^{-1} + \cdots)
\\[11pt]
&=\ 
s_x + 2^{-(x + 1)}(2)
\\[11pt]
&=\ 
s_x + 2^{-x}.
\end{align*}
Now quote $\S13$, $\#18$.
\end{x}
\vspace{0.3cm}

\begin{x}{\small\bf REMARK} \ %09
Assume that \mA is not recursive $-$then for certain \mA, it is possible to find an injective $f \in \fR$ such that \mA equals the range of $f$, 
hence
\[
\alpha_A 
\ = \ 
\sum\limits_{n \in \N} \hsx 2^{-f(n)}.
\]
Put
\[
q_x 
\ = \ 
\sum\limits_{n \leq x} \hsx 2^{-f(n)} \qquad (x = 0, 1, 2, \ldots).
\]
Then $q_0$, $q_1$, $q_2, \ldots$ is a bounded increasing sequence of rational numbers whose limit $\alpha_A$ is not computable.
\end{x}
\vspace{0.3cm}

We have worked thus far with 
\[
\R_\fR \ \cap \ [0,\infty[.
\]
Replace now $\R_\fR$ by $\R_{\fP\fR}$ $-$then it is still the case that
\[
\R_{\fP \fR}^b \ \subset \ \R_{\fP \fR} \hsx \cap \hsx [0,\infty[
\]
but it is no longer true that
\[
\R_{\fP \fR}^b \ = \ \R_{\fP \fR} \hsx \cap \hsx [0,\infty[\hsx .
\]
\\
%%----------------------------------------------------------------------------------------------06

\begin{x}{\small\bf THEOREM} \ %10
$\R_{\fP \fR}$ is a field (cf. $\S 8$, $\# 15$).
\end{x}
\vspace{0.3cm}

\begin{x}{\small\bf LEMMA}%11
\footnote[2]{Qingliang Chen et al. \textit{Mathematical Logic Quarterly} {\bf{53}} (2007), pp. 365-380.}
\ 
For each $b > 1$, there are $x$, $y \in  \ \R_{\fP \fR}^b$ such that $x + y \not\in \ \R_{\fP \fR}^b$ .
\end{x}
\vspace{0.3cm}

Consequently the containment
\[
\R_{\fP \fR}^b \ \subset \ \R_{\fP \fR} \cap \hsx [0,\infty[
\]
is proper.
\\[-0.25cm]

One can also consider the relationship among the $\R_{\fP \fR}^b$ for different $b$.
\\[-0.25cm]

\begin{x}{\small\bf THEOREM} \ %12
Let $b$, $d > 1$ $-$then
\[
\R_{\fP \fR}^b \ \subset \ \R_{\fP \fR}^d
\]
iff $d$ divides a power of $b$, i.e., iff there exist $k$, $s \in \ \N$ such that $b^k = s d$.  
\\[-0.5cm]

[For details and references, see the paper of Qingliang Chen et al.]
\end{x}

%%%%%%%%%%%%%%%%%%%%%%%%%%%%%%%%%%%%%%
%%%%%%%%%%%%%%%%%%%%%%%%%%%%%%%%%%%%%%
%%%%%%%%%%%%%%%%%%%%%%%%%%%%%%%%%%%%%%

%12
%\include{__refs}
%\newpage
\setcounter{page}{1}
\renewcommand{\thepage}{Index-\arabic{page}}
\backmatter
\bibliography{}
\printindex
\end{document}